\documentclass[final]{siamltex}

\usepackage{amsmath}
\usepackage{enumitem}
\usepackage{amssymb, latexsym, mathrsfs}
\usepackage{fancybox, graphicx, subfigure}
\usepackage{endnotes}
\usepackage[usenames, dvipsnames]{color}
\usepackage[colorlinks=true,
        raiselinks=true,
        pagecolor=magenta,
        linkcolor=myblue,
        citecolor=mygreen,
        urlcolor=mygreen,
        pdfauthor=Debasish Chatterjee,
        pdftitle={},
        pdfkeywords={},
        pdfsubject={Technical Report},
        plainpages=false]{hyperref}

\definecolor{myred}{rgb}{0.6, 0, 0}
\definecolor{mygreen}{rgb}{0, 0.5, 0}
\definecolor{myblue}{rgb}{0, 0, 0.5}
\definecolor{mycyan}{rgb}{0, 0.5, 0.5}

\newcommand{\R}{\ensuremath{\mathbb{R}}}
\newcommand{\N}{\ensuremath{\mathbb{N}}}
\newcommand{\Nz}{\ensuremath{\mathbb{N}_0}}
\newcommand{\posR}{\ensuremath{\R_{\ge 0}}}

\newcommand{\PSet}{\ensuremath{\mathcal{P}}}

\newcommand{\CSet}{\ensuremath{\mathcal{U}}}

\newcommand{\ra}{\ensuremath{\rightarrow}}
\newcommand{\Ra}{\ensuremath{\;\Longrightarrow\;}}
\newcommand{\lra}{\ensuremath{\longrightarrow}}
\newcommand{\eps}{\ensuremath{\varepsilon}}
\newcommand{\fa}{\ensuremath{\forall\,}}
\renewcommand{\le}{\ensuremath{\leqslant}}
\renewcommand{\ge}{\ensuremath{\geqslant}}
\renewcommand{\mapsto}{\ensuremath{\longmapsto}}

\newcommand{\therex}{\ensuremath{\exists\,}}
\newcommand{\setmin}{\ensuremath{\!\smallsetminus\!}}

\newcommand{\ClassK}{\ensuremath{\mathcal{K}}}
\newcommand{\ClassKinfty}{\ensuremath{\mathcal{K}_{\infty}}}
\newcommand{\ClassKL}{\ensuremath{\mathcal{KL}}}

\newcommand{\epower}[1]{\ensuremath{\mathrm{e}^{#1}}}
\newcommand{\norm}[1]{\ensuremath{\left\lVert #1 \right\rVert}}

\newcommand{\Expec}[1]{\ensuremath{\mathsf{E}\!\left[\vphantom{\big|}#1\vphantom{\big|}\right]}}
\newcommand{\CExpec}[2]{\ensuremath{\mathsf{E}^{#2}_{\vphantom{T}}\!\!\left[\vphantom{\big|}#1\vphantom{\big|}\right]}}
\newcommand{\CProb}[2]{\ensuremath{\mathsf{P}^{#2}_{\vphantom{T}}\!\!\left(\vphantom{\big|}#1\vphantom{\big|}\right)}}
\newcommand{\Prob}[1]{\ensuremath{\mathsf{P}\!\left(\vphantom{\big|}#1\vphantom{\big|}\right)}}
\newcommand{\LieD}[2]{\ensuremath{\mathrm L_{#1}{#2}}}
\newcommand{\indic}[1]{\ensuremath{\boldsymbol{1}_{#1}}}
\newcommand{\abs}[1]{\ensuremath{\left\lvert{#1}\right\rvert}}

\newcommand{\mrm}[1]{\ensuremath{\mathrm{#1}}}

\newcommand{\Lp}[1]{\ensuremath{\boldsymbol{L}_{#1}}}

\newcommand{\secref}[1]{\S\ref{#1}}

\newcommand{\transp}{\ensuremath{^{\scriptscriptstyle{\mathrm T}}}}
\newcommand{\sigalg}{\ensuremath{\mathfrak{F}}}
\newcommand{\cardP}{\ensuremath{\mathrm{N}}}

\newcommand{\xz}{\ensuremath{x_0}}
\newcommand{\sigmaz}{\ensuremath{\sigma_0}}

\newcommand{\drv}{\ensuremath{\,\mathrm{d}}}
\newcommand{\ol}{\overline}

\newcommand{\wt}{\widetilde}

\newcommand{\gasas}{{\sc gas}~a.s.}
\newcommand{\gasm}{{\sc gas-m}}

\newcommand{\gasp}{{\sc gas-p}}

\newcommand{\cadlag}{c{\`a}dl{\`a}g}
\renewcommand{\subset}{\ensuremath{\subseteq}}

\newcommand{\mx}{\ensuremath{\vee}}
\newcommand{\mn}{\ensuremath{\wedge}}

\newcommand{\RemarkEnd}{\hspace{\stretch{1}}{$\vartriangleleft$}}

\newtheorem{defn}[theorem]{Definition}
\newtheorem{remark}[theorem]{Remark}

\newtheorem{assumption}[theorem]{Assumption}

\newcommand{\iss}{\textsc{iss}}

\title{Stabilizing Randomly Switched Systems\thanks{This work was supported by NSF CSR program (Embedded \& Hybrid systems area) under grant NSF CNS-0614993.}}
\author{Debasish Chatterjee\thanks{Automatic Control Laboratory, ETL K24, ETH Z\"urich, 8092 Z\"urich, Switzerland (Email: {\tt chatterjee@control.ee.ethz.ch}).} \and Daniel Liberzon\thanks{144 Coordinated Science Laboratory, University of Illinois at Urbana-Champaign, Urbana, IL 61801, USA (Email: {\tt liberzon@uiuc.edu}).}}

\begin{document}

	\maketitle

	\begin{abstract}
		This article is concerned with stability analysis and stabilization of randomly switched systems under a class of switching signals. The switching signal is modeled as a jump stochastic (not necessarily Markovian) process independent of the system state; it selects, at each instant of time, the active subsystem from a family of systems. Sufficient conditions for stochastic stability (almost sure, in the mean, and in probability) of the switched system are established when the subsystems do not possess control inputs, and not every subsystem is required to be stable. These conditions are employed to design stabilizing feedback controllers when the subsystems are affine in control. The analysis is carried out with the aid of multiple Lyapunov-like functions, and the analysis results together with universal formulae for feedback stabilization of nonlinear systems constitute our primary tools for control design.
	\end{abstract}



	\pagestyle{myheadings}
	\thispagestyle{plain}
	\markboth{D.~CHATTERJEE {\sc and} D.~LIBERZON}{STABILIZING RANDOMLY SWITCHED SYSTEMS}

	\section{Introduction}
		\label{s:intro}
		A randomly switched system has two ingredients, namely, a family of subsystems and a random switching signal. In this article we are interested in finding conditions for stochastic stability of randomly switched systems. Our approach consists of identifying key properties of the family of subsystems and the switching signal, and finding conditions to connect them such that the switched system has the desired characteristics. We concentrate on stability almost surely and in expectation. Since each of these implies stability in probability~\cite{ref:khasminskiiStochastic, ref:kushnerStocStab, ref:kozinsurvey}, our results immediately provide sufficient conditions for weak stability in probability of the systems under consideration; we also demonstrate that the conditions are sufficient for strong stability in probability.

		The basic structure of our main analysis results is as follows. The first step involves extracting properties which quantitatively express stability characteristics of the subsystems. This is carried out with the help of multiple Lyapunov functions. The method of multiple Lyapunov functions was developed originally in the context of deterministic switched systems, and is discussed in detail in, e.g.,~\cite[Chapter~3]{ref:liberzonbk}. This method is effective in quantitatively capturing the degree of stability (or instability) of the subsystems. The second step involves extracting key properties of the switching signal. These properties are variously captured by the probability mass function of its rate of switching, the probability distribution of its jump destinations, distribution of holding times between switching instants, etc. Finally, the characteristics of the switched system generated by the switching signal from the family of subsystems are captured by inequalities which connect the above two sets of properties.

		Research on randomly switched systems has concentrated mostly on the case of Markovian switching signals---the discrete state evolves according to a continuous-time Markov chain; see e.g.,~\cite{ref:davisMarkovModelsOptbk, ref:maoasystabdist} and the references cited therein. The central idea behind arriving at stability conditions revolves around employing the generator of the Markov process and extracting certain nonnegative supermartingales that converge to zero in expectation. This method in turn is based on the martingale problem~\cite[Chapter~5]{ref:ethierkurtzMP} corresponding to the Markov process. In its simplest form, if $(X_t)_{t\ge 0}$ is the underlying Markov process with generator $\mathcal L$, ($X$ consists of both the continuous and the discrete states,) then for every measurable and bounded real-valued function $V$ on the state space, the process $(Y_t)_{t\ge 0}$ defined by $Y_t := V(X_t) - \int_0^t \mathcal LV(X_s)\mrm ds$ is a martingale. A pointwise inequality which bounds $\mathcal LV(\cdot)$ on the state space may be imposed, and with the help of this one can draw conclusions about stability properties of the system by analyzing the martingale above; this analysis becomes particularly simple if $V$ is a nonnegative Lyapunov-like function.

		The martingale approach described above can be applied to switched systems in which the switching signals are general point processes with intensity functions satisfying certain standard measurability conditions; see e.g.,~\cite{ref:bremaudPointProc} for further details on the measurability conditions. These intensity functions appear in the expression of $\mathcal LV$ in place of the usual Markov transition intensity matrix, and hereafter the analysis follows that of the Markovian case. However, for non-Markovian switching signals, it is not easy to employ this technique; for instance, if the holding times between consecutive switching instants are independent and identical uniform random variables, obtaining expressions of these intensity functions is difficult. The methods we propose here apply to semi-Markovian switching signals, do not depend on martingale analysis, and arrive at the results directly by employing what we think are less involved and more intuitively appealing techniques. Existing work on stability of stochastic switched systems includes~\cite{ref:maoSwitchingbk, ref:yuanlygerosSDEMarksw, ref:ugrinovskiiRA, ref:lygerosReachQuestions, ref:boukasSSS, ref:davisMarkovModelsOptbk, ref:fengMarkovJump92, ref:jichizeck90, ref:hespStocHybCommNet, ref:battilottiDwellTimeMJS}; see also~\cite[Chapter~1]{ref:myphdthesis} for a survey techniques employed in these articles.

		Analysis results obtained via our approach, including those reported in our earlier article~\cite{ref:ransw} where each subsystem was required to be stable, have conceptual analogs in deterministic switched systems theory. The approach pursued in~\cite{ref:ransw} and in the current article is derived from the method of multiple Lyapunov functions developed in the context of deterministic switched systems, see e.g.,~\cite[Chapter~3]{ref:liberzonbk} for an extensive discussion. Stability of individual subsystems and a slow switching condition are the important features of these deterministic results. In this article our results involving unstable subsystems employ certain probabilistic characteristics of the switching signal in addition to slow switching; their conceptual analogs in deterministic switched systems literature are comparatively less known, with the exception of~\cite{ref:michelunstabswsys}.

		With our analysis results in hand, we turn to control synthesis and derive explicit controller formulas which ensure stability of the switched system in closed loop. In this context, there naturally arise two distinct cases: one in which the controller has full knowledge of the switching signal at each instant of time, and the other in which the controller is totally unaware of the switching signal. We examine the distinctive features of each of these two cases and propose control synthesis strategies by employing \emph{universal formulae}~\cite{ref:sontagunivformula, ref:linunivformulabddcon, ref:linunivformularestrcon, ref:malisoffunivformulaMinkball} for nonlinear feedback stabilization. The advantages of our approach are evident here, for one does not need to design a controller separately for the switched system if there already exist \emph{control-Lyapunov functions} for each individual subsystem; then, off-the-shelf controllers employing universal formulae are easily designed, and a modular organization of the controller synthesis stage is facilitated.

		The article unfolds as follows. \secref{s:prelims} presents the system model with no inputs and the stability concepts under consideration. The main analysis results appear in~\S\S\ref{s:mainres}, \ref{s:general}, and~\ref{s:gasp}, and their proofs are given in~\secref{s:proofs}. The controller synthesis results are presented in~\secref{s:syn}. We conclude in~\secref{s:concl} with a brief discussion of possible channels of further investigation.

		Some notation: Let $\posR$ denote the nonnegative half-line $[0, \infty[$, $\N = \{1, 2, \ldots\}$, $\Nz := \N\cup\{0\}$, and let $\norm{\cdot}$ denote the Euclidean norm.

	\section{Preliminaries}
		\label{s:prelims}
		We define the family of systems
        \begin{equation}
            \dot x = f_i(x),\qquad i\in\PSet,
            \label{e:ssysfam}
        \end{equation}
        where the state $x\in\R^n$, $\PSet$ is a finite index set of $\cardP$ elements: $\PSet = \{1, \ldots, \cardP\}$, the vector fields $f_i:\R^n\lra\R^n$ are locally Lipschitz, and $f_i(0) = 0$, $i\in\PSet$.

		Let $(\Omega, \sigalg, {\mathsf P})$ be a complete probability space. Let $\sigma := (\sigma(t))_{t\ge 0}$ be a \cadlag{} (i.e., right-continuous and possessing limits from the left) stochastic process taking values in $\PSet$, with $\sigma(0)$ completely known. The process $\sigma$ is by definition measurable~\cite[Chapter~1]{ref:revuzyorCMBM}. Let the discontinuity points of $\sigma$ be denoted by $\tau_i, \;i\in\N$, and let $\tau_0 := 0$ by convention. The filtration $(\sigalg_t)_{t\ge 0}$ generated by $\sigma$ is right-continuous~\cite[Theorem~T26, p.~304]{ref:bremaudPointProc}, and we augment $\sigalg_0$ with all $\mathsf P$-null sets. As a consequence of the hypotheses of our results, the sequence $(\tau_i)_{i\in\Nz}$ is almost surely divergent, i.e., $\sigma$ is nonexplosive. The {\em randomly switched system} generated by this {\em switching signal} $\sigma$ from the family~\eqref{e:ssysfam} is
        \begin{equation}
            \dot x = f_\sigma(x), \qquad (x(0), \sigma(0)) = (\xz, \sigma_0), \quad t\ge 0.
            \label{e:ssysdef}
        \end{equation}
		We assume that there are no jumps in the state $x$ at the points of discontinuity of the switching signal; we shall henceforth refer to these points as the switching instants. The above hypotheses on the system~\eqref{e:ssysdef} and $\sigma$ ensure that standard conditions for the existence and uniqueness of an absolutely continuous solution in the sense of Carath{\'e}odory~\cite{ref:filippovbk}, over a nontrivial time interval containing $0$, are fulfilled for almost every sample path. Existence and uniqueness of a global solution will follow from the hypotheses of our results. We let $x(\cdot)$ denote this solution. For $\xz = 0$, the solution to~\eqref{e:ssysdef} is identically $0$ for every $\sigma$; we shall ignore this trivial case in the sequel. Standard arguments (see e.g,~\cite[Chapter~1]{ref:myphdthesis}) show that the solution process $x(\cdot)$ of~\eqref{e:ssysdef} is an $(\sigalg_t)_{t\ge 0}$-adapted process.

		Recall~\cite{ref:borovkovbk} that for $\lambda > 0$, an exponential-$(\lambda)$ random variable $\xi$ has the distribution function $\Prob{\xi \le s} = 1-\epower{-\lambda s}$ for $s \ge 0$, and $0$ otherwise; for $T > 0$, a uniform-$(T)$ random variable $\xi$ has the distribution function $\Prob{\xi \le s} = 0$ if $s < 0$, $s/T$ if $s\in[0, T]$, and $1$ otherwise. A continuous function $\alpha:\posR\lra\posR$ is of class-$\mathcal K$ (we write $\alpha\in\ClassK$) if it vanishes at $0$ and is monotone strictly increasing. A continuous function $\beta:\posR\times\posR\lra\posR$ is of class-$\ClassKL$ (we write $\beta\in\ClassKL$) if $\beta(r, \cdot)$ is monotone strictly decreasing for each fixed $r$, and $\beta(\cdot, s)$ is of class-$\ClassK$ for each fixed $s$; we write $\beta\in\ClassKL$.

		We focus on the following two properties of~\eqref{e:ssysdef}; see e.g.,~\cite{ref:khasminskiiStochastic}.

		\begin{defn}
			The system~\eqref{e:ssysdef} is said to be \emph{globally asymptotically stable almost surely} ({\sc gas}~a.s.) if the following two properties are simultaneously verified:
			\begin{enumerate}[label={\rm (AS\arabic*)}, align=left, leftmargin=*]
				\item $\Prob{\fa\eps > 0\;\;\therex\delta > 0 \text{ such that }\displaystyle{\norm{\xz} < \delta\Ra \sup_{t\ge 0}\norm{x(t)} < \eps}} = 1$;
				\item $\Prob{\fa r, \eps' > 0\;\;\therex T \ge 0\text{ such that }\displaystyle{\norm{\xz} < r\Ra \sup_{t\ge T}\norm{x(t)} < \eps'}} = 1$.
			\end{enumerate}
			\label{d:gasas}
		\end{defn}

		Let us note that this property is well-defined because each of the sets appearing inside the measure $\mathsf P$ is $\sigalg$-measurable due to continuity of $x(\cdot)$.

		\begin{defn}
			The system~\eqref{e:ssysdef} is said to be \emph{$\alpha$-globally asymptotically stable in the mean} ($\alpha$-\gasm{}) for a function $\alpha\in\ClassK$ if the following two properties are simultaneously verified:
			\begin{enumerate}[label={\rm (SM\arabic*)}, align=left, leftmargin=*]
				\item $\fa\eps > 0\;\;\therex\wt\delta > 0$ such that $\displaystyle{\norm{\xz} < \wt\delta\Ra\sup_{t\ge 0}\Expec{\alpha(\norm{x(t)})} < \eps}$;\medskip
				\item $\fa r, \eps' > 0\;\; \therex \wt T \ge 0$ such that $\displaystyle{\norm{\xz} < r\Ra\sup_{t\ge \wt T}\Expec{\alpha(\norm{x(t)})} < \eps'}$.
			\end{enumerate}
			\label{d:gasm}
		\end{defn}

		Stability definitions in deterministic systems literature usually involve just the norm of the state. The presence of the function $\alpha$ in Definition~\ref{d:gasm} allows some measure of flexibility in the sense that one need not worry about bounds for just the expectation of the norm of the state, i.e., $\Lp 1$-stability. Frequently one employs Lyapunov functions which are polynomial functions of the states, and with the aid of conditions such as (V1) in Assumption~\ref{a:V} below, stronger bounds in terms of the $\Lp p$ ($p > 1$) norms of the state are obtained. For instance, quadratic Lyapunov functions yield bounds for mean-square or $\Lp 2$-stability, which is stronger than $\Lp 1$-stability.

		Our analysis results employ a family of Lyapunov functions, one for each subsystem.
		The following assumption collects the properties we shall require from the members of this family of Lyapunov functions.\footnote{Strictly speaking we should call them ``Lyapunov-like functions,'' because their gradients do not necessarily decrease along the corresponding system trajectories. For simplicity we shall adhere to the term ``Lyapunov functions'' in the sequel.} For notational brevity, we let $\LieD{f}{V}(x)$ denote the Lie derivative of a differentiable function $V:\R^n\lra\R$ along a vector field $f:\R^n\lra\R^n$, i.e., $\LieD{f}V(x) := \left\langle\nabla_x V(x)\vphantom{\big|}, \vphantom{\big|}f(x)\right\rangle$.

		\begin{assumption}
			\label{a:V}
			There exist a family of continuously differentiable real-valued functions $\{V_i\}_{i\in\PSet}$ on $\R^n$, functions $\alpha_1, \alpha_2\in\ClassKinfty$, numbers $\mu > 1$ and $\lambda_i\in\R$, $i\in\PSet$, such that for all $x\in\R^n$ and $i, j\in\PSet$,
			\begin{enumerate}[label={\rm (V\arabic*)}, align=left, leftmargin=*]
				\item $\alpha_1(\norm x)\le V_i(x)\le \alpha_2(\norm x)$;
				\item $\displaystyle{\LieD{f_i}{V_i}(x)\le -\lambda_{i} V_i(x)}$;
				\item $V_i(x) \le \mu V_j(x)$.
			\end{enumerate}
		\end{assumption}

		\begin{remark}
			{\rm
				(V1) is a fairly standard hypothesis, ensuring each $V_i$ is positive definite and radially unbounded. The condition in (V2) keeps track of the growth of $i$-th Lyapunov function $V_i$ along the vector field $f_i$ of the $i$-th subsystem; the parameter $\lambda_{i}$ provides a quantitative estimate of this growth rate. The right-hand side of the inequality in (V2) being a linear function of $V_i$ is no loss of generality, see e.g.,~\cite[Theorem~2.6.10]{ref:lakshmikanthamDifferentialIntegralInequalities1} for details. (V3) certainly restricts the class of functions that the family $\{V_i\}_{i\in\PSet}$ can belong to; however, this hypothesis is commonly employed in the deterministic context~\cite[Chapter~3]{ref:liberzonbk}. Quadratic Lyapunov functions universally utilized in the case of linear subsystems always satisfy this hypothesis.
			}
		\end{remark}

	\section{Main Results}
		\label{s:mainres}
		In this section we present our main results providing sufficient conditions for {\sc gas}~a.s.\ and $\alpha_1$-\gasm{} of randomly switched systems under two different classes of switching signals. The switching signals described here are fairly general and are quite natural to consider.

		We let $(S_i)_{i\in\N}$, $S_i := \tau_i - \tau_{i-1}$ be the sequence of holding times, where $(\tau_i)_{i\in\N}$ is the sequence of discontinuity points of $\sigma$.

		\begin{defn}\mbox{}
			\label{d:sigmaclasses}
			We say that the switching signal $\sigma$ belongs to
			\begin{itemize}[label=\textbullet, leftmargin=*]
				\item \emph{class EH} if:
					\begin{enumerate}[label={\rm (EH\arabic*)}, align=left, leftmargin=*]
						\item the sequence $(S_i)_{i\in\N}$ of holding times is a collection of independent and identically distributed (i.i.d) random variables, with $S_i$ an exponential-$(\lambda)$ random variable, $\lambda > 0$;
						\item $\therex q_i\in[0, 1]$, $i\in\PSet$, such that $\fa j\in\N$, $\Prob{\sigma(\tau_j) = i\big|(\sigma(\tau_k))_{k=0}^{j-1}} = q_i$;
						\item the sequences $(S_i)_{i\in\N}$ and $(\sigma(\tau_i))_{i\in\Nz}$ are mutually independent.
					\end{enumerate}
					\medskip
				\item \emph{class UH} if:
					\begin{enumerate}[label={\rm (UH\arabic*)}, align=left, leftmargin=*]
						\item the sequence $(S_i)_{i\in\N}$ of holding times is a collection of i.i.d random variables, with $S_i$ a uniform-$(T)$ random variable, $T > 0$;
						\item $\therex q_i\in[0, 1]$, $i\in\PSet$, such that $\fa j\in\N$, $\Prob{\sigma(\tau_j) = i\big|(\sigma(\tau_k))_{k=0}^{j-1}} = q_i$;
						\item the sequences $(S_i)_{i\in\N}$ and $(\sigma(\tau_i))_{i\in\Nz}$ are mutually independent.
					\end{enumerate}
			\end{itemize}
		\end{defn}
		\medskip

		The following are our main results; their proofs are provided in~\secref{s:proofs}.
		\medskip

		\begin{theorem}
			\label{t:gasaspoisson}
			The system~\eqref{e:ssysdef} is {\sc gas}~a.s.\ if
			\begin{enumerate}[label={\rm (E\arabic*)}, align=left, leftmargin=*]
				\item Assumption~\ref{a:V} holds;
				\item the switching signal $\sigma$ belongs to class EH as defined in Definition~\ref{d:sigmaclasses};
				\item $\lambda_i + \lambda > 0\quad \fa i\in\PSet$;
				\item $\displaystyle{\sum_{i\in\PSet}\left(\frac{\mu q_i}{1+\lambda_i/\lambda}\right) < 1}$.
			\end{enumerate}
		\end{theorem}

		\begin{corollary}
			\label{c:gasmpoisson}
			The system~\eqref{e:ssysdef} is $\alpha_1$-\gasm{} under the hypotheses of Theorem~\ref{t:gasaspoisson}.
		\end{corollary}

		\begin{theorem}
			\label{t:gasasunif}
			The system~\eqref{e:ssysdef} is {\sc gas}~a.s.\ if
			\begin{enumerate}[label={\rm (U\arabic*)}, align=left, leftmargin=*]
				\item Assumption~\ref{a:V} holds;
				\item the switching signal $\sigma$ belongs to class UH as defined in Definition~\ref{d:sigmaclasses};
				\item $\displaystyle{\sum_{i\in\PSet}\left(\frac{\mu q_i\left(1-\epower{-\lambda_i T}\right)}{\lambda_i T}\right) < 1}$.
			\end{enumerate}
		\end{theorem}

		\begin{corollary}
			\label{c:gasmunif}
			The system~\eqref{e:ssysdef} is $\alpha_1$-\gasm{} under the hypotheses of Theorem~\ref{t:gasasunif}.
		\end{corollary}

		\medskip

		\begin{remark}
			\label{r:nozeno}
			{\rm 
				Let us first note that switching signals of class EH and UH are nonexplosive, i.e., there are finitely many jumps on finite-length intervals of time almost surely. Indeed, it follows immediately from the Strong Law of Large Numbers~\cite[Theorem~7, p.~64]{ref:raoProbTheo} that since $(S_i)_{i\in\N}$ is i.i.d and $\Expec{S_i}\in\;]0, \infty[$ for switching signals belonging to either class EH or UH, almost surely the $\nu$-th jump instant $\tau_\nu = \sum_{i=1}^\nu S_i \ra\infty$ as $\nu\ra\infty$. It is also readily seen that switching cannot stop after a finite time, for then $S_j = \infty$ for some $j$, and the probability of the event $\{S_j = \infty\text{ for some $j$}\}$ is $0$.
			}
		\end{remark}

		\begin{remark}
			\label{r:theorempoisson}
			{\rm 
				Let us examine the statement of Theorem~\ref{t:gasaspoisson} in some detail. Firstly, note that by (E1) not all subsystems are required to be stable, i.e., for some $i\in\PSet$, $\lambda_i$ can be negative; then (V2) provides a measure of the rate of instability of the corresponding subsystems. Secondly, note that condition (E3) is always satisfied if each $\lambda_i > 0$. However, if $\lambda_i < 0$ for some $i\in\PSet$, then (E3) furnishes a maximum instability margin of the corresponding subsystems that can still lead to {\sc gas}~a.s.\ of~\eqref{e:ssysdef}. Intuitively, in the latter case, the process $N_\sigma(t, 0)$ must switch fast enough (which corresponds to $\lambda > 0$ being large enough,) so that the unstable subsystems are not active for too long. Potentially this fast switching may have a destabilizing effect. Indeed, it may so happen that for a given $\mu$, a fixed probability distribution $\{q_i\}_{i\in\PSet}$, and a choice of functions $\{V_i\}_{i\in\PSet}$, (E3) and (E4) may be impossible to satisfy simultaneously, due to a very high degree of instability of even one subsystem for which the corresponding $q_i$ is also large. Then we need to search for a different family of functions $\{V_i\}_{i\in\PSet}$ for which the hypotheses hold. Thirdly, (E4) connects the properties of deterministic subsystem dynamics, furnished by the family of Lyapunov functions satisfying Assumption~\ref{a:V}, with the properties of the stochastic switching signal. From (E4) it is clear that larger degrees of instability of a subsystem (small $\lambda_i$) can be compensated by a smaller probability of the switching signal activating the corresponding subsystem.
			}
		\end{remark}

		\begin{remark}
			\label{r:theoremunif}
			{\rm 
				Let us make some observations about the statement of Theorem~\ref{t:gasasunif}. Once again, just like Theorem~\ref{t:gasaspoisson}, note that by (U1) not all subsystems are required to be stable, i.e., for some $i\in\PSet$, $\lambda_i$ can be negative. (U3) connects the properties of deterministic subsystem dynamics, furnished by the family of Lyapunov functions satisfying Assumption~\ref{a:V}, with the properties of the stochastic switching signal. Also from (U3) it is clear that larger degrees of instability (larger $\lambda_i$) of a subsystem can be compensated by a smaller probability (smaller $q_i$) of the switching signal activating the corresponding subsystem. Notice that a switching signal of class UH is semi-Markov~\cite[Section~20.4]{ref:borovkovbk}. There is a strong dependence on past history due to the uniform holding times. Indeed, at an arbitrary instant of time $t$ we need to know how long ago the last jump occurred in order to compute the probability distribution of the next jump instant after $t$.
			}
		\end{remark}

		\begin{remark}
			\label{r:comppoissonunif}
			{\rm
				It may be observed that Theorem~\ref{t:gasaspoisson} requires a larger set of hypotheses compared to Theorem~\ref{t:gasasunif}; however, this is only natural. Indeed, the switching signal in the latter case is constrained to switch at least once in $T$ units of time, whereas no such hard constraint is present on the switching signal in the former case. We observed in Remark~\ref{r:theorempoisson} that it is necessary for the switching signal to switch fast enough if there are unstable subsystems in the family~\eqref{e:ssysfam}, which necessitated the condition (E3). This fast switching is automatic if $\sigma$ is of class UH, provided $T$ is related to the instability margin of the subsystems in a particular way. The condition (U3) captures this relationship, for, observe that if $\lambda_i$ is negative and large in magnitude for some $i\in\PSet$, the ratio $\left(1-\epower{-\lambda_i T}\right)/(\lambda_i T)$ is smaller for smaller $T$, and a smaller ratio is better for {\sc gas}~a.s.\ of~\eqref{e:ssysdef}. Also for a given $T$, large and positive $\lambda_i$'s (i.e., subsystems with high margins of stability) make the aforesaid ratio small.
			}
		\end{remark}

	\section{A Generalization}
		\label{s:general}
		The results in~\secref{s:mainres} fall short of being completely satisfactory. In particular, the assumption of the jump destinations process $(\sigma(\tau_i))_{i\in\N}$ being memoryless (assumptions (EH2) and (UH2)) is perhaps the most restrictive. As we observed in Remark~\ref{r:theoremunif}, switching signals of class UH fall in the class of semi-Markov processes, in fact trivially so, due to the memoryless nature of the discrete jump-destination process $(\sigma(\tau_i))_{i\in\N}$. However, it would be better if we could handle the Markovian jump destination case by keeping the other two hypotheses intact. In this section we do that, namely, include those switching signals for which the process $(\sigma(\tau_i))_{i\in\N}$ is a discrete-time Markov chain. Although the results given in this section are not the most general possible, they are intended to highlight the directions of possible generalizations that can be made in our framework.

		\begin{assumption}
			\label{a:Vadd}
			There exist a family of continuously differentiable real-valued functions $\{V_i\}_{i\in\PSet}$ on $\R^n$, functions $\alpha_1, \alpha_2\in\ClassKinfty$, numbers $\mu > 1$ and $\lambda_{i, j}\in\R$, $i, j\in\PSet$, such that for all $x\in\R^n$ and $i, j\in\PSet$,
			\begin{enumerate}[label={\rm (V\arabic*$'$)}, align=left, leftmargin=*]
				\item {\rm (V1)} of Assumption~\ref{a:V} holds;
				\item $\displaystyle{\LieD{f_j}{V_i}(x)\le -\lambda_{i, j} V_i(x)}$;
				\item {\rm (V3)} of Assumption~\ref{a:V} holds.
			\end{enumerate}
		\end{assumption}
		\medskip

		\begin{defn}
			\label{d:sigmaclassgen}
			We say that the switching signal $\sigma$ belongs to \emph{class GH} if:
			\begin{enumerate}[label={\rm (GH\arabic*)}, align=left, leftmargin=*]
				\item the sequence $(S_i)_{i\in\N}$ of holding times is an i.i.d collection of random variables, with $\Expec{S_i} < \infty$;
				\item the process $(\sigma(\tau_i))_{i\in\Nz}$ is a discrete-time Markov chain with initial probability vector\footnote{Here $\delta_{\{j\}}$ denotes the Dirac measure concentrated on $\{j\}$.} $\delta_{\{\sigma_0\}}$ and transition probability matrix $P = [p_{i, j}]_{\PSet\times\PSet}^{\vphantom T}$;
				\item $(S_i)_{i\in\N}$ is independent of $(\sigma(\tau_i))_{i\in\Nz}$.
			\end{enumerate}
		\end{defn}

		Switching signals belonging to class GH are semi-Markov~\cite[Section~20.4]{ref:borovkovbk}. In the most general case of a semi-Markov process, the sequence $(S_i)_{i\in\N}$ in (GH1) may be such that the distribution of $S_i$ depends on both $\sigma(\tau_{i-1})$ and $\sigma(\tau_i)$, $i\in\N$. Our objective here is to illustrate some new techniques, and hence we shall retain the simpler condition (GH1) at the expense of lesser generality. The condition (GH2) imposes a discrete-time Markovian structure on the process $(\sigma(\tau_i))_{i\in\Nz}$, and the condition (GH3), though not the most general, is a standard hypothesis for semi-Markov processes.

		\begin{theorem}
			\label{t:gasasgen}
			The system~\eqref{e:ssysdef} is {\sc gas}~a.s.\ if
			\begin{enumerate}[label={\rm (G\arabic*)}, align=left, leftmargin=*]
				\item Assumption~\ref{a:Vadd} holds;
				\item the switching signal $\sigma$ belongs to class GH as defined in Definition~\ref{d:sigmaclassgen};
				\item $\therex \theta\in[0, 1[$ such that
					\[
						\max_{i\in\PSet}\sum_{j\in\PSet}\left(\vphantom{\sum}\mu p_{i, j}\Expec{\epower{-\lambda_{j, i}S_k}}\right) \le \theta.
					\]
			\end{enumerate}
		\end{theorem}

		\begin{remark}
			{\rm 
				Switching signals of class GH are nonexplosive, and switching cannot stop in finite time, as can be seen by following the same line of reasoning as in Remark~\ref{r:nozeno}.
			}
		\end{remark}

		\begin{remark}
			{\rm 
				Note that Theorem~\ref{t:gasasgen} is conceptually quite different from the results of~\secref{s:mainres}. Indeed, the condition (G3) involves the growth rate of a Lyapunov function along every subsystem, in contrast to the results in~\secref{s:mainres}, where we only kept track of the growth rate of each Lyapunov function along the trajectories of the corresponding subsystem. This additional factor is due to the Markovian nature of the jump destination process $(\sigma(\tau_i))_{i\in\Nz}$, and quite naturally the transition probabilities $p_{i, j}, i, j\in\PSet$ appear in (G3). Also, the condition (V2$'$) requires us to keep track of the behavior of every Lyapunov function at once; in a way we quantify how each subsystem relates to the others through the inequality in (V2$'$). This is a deviation from our philosophy of decoupling the properties of the switching signal from the properties of the individual subsystems at first and then connecting them. The Markovian nature of the jump destination process in Theorem~\ref{t:gasasgen} does not seem to entirely allow this separation.
			}
		\end{remark}

	\section{An Excursion into Global Asymptotic Stability in Probability}
	\label{s:gasp}
		Among the several notions of stochastic stability in the literature, one particular notion that encodes uniform behavior of system trajectories is strong global asymptotic stability in probability (s-\gasp{}). Recall~\cite{ref:khasminskiiStochastic} that
		\begin{defn}
			\label{d:gasp}
			The system~\eqref{e:ssysdef} is \emph{strongly globally asymptotically stable in probability} if the following two properties are simultaneously verified:
			\begin{enumerate}[label={\rm (\roman*)}, align=right, leftmargin=*, widest=ii]
				\item $\fa \eta\in\;]0, 1[\;\; \fa \eps > 0\;\therex\delta > 0$ such that $\displaystyle{\norm{\xz} < \delta\Ra\Prob{\sup_{t\ge 0}\norm{x(t)} > \eps} \le \eta}$;
				\item $\fa \eta'\in\;]0, 1[\;\;\fa r, \eps' > 0\;\therex T > 0$ such that $\displaystyle{\norm{\xz} < r\Ra\Prob{\sup_{t\ge T}\norm{x(t)} > \eps'} \le \eta'}$.
			\end{enumerate}
		\end{defn}

		Let us note that each of the sets inside the measure $\mathsf P$ in (i) and (ii) above is $\sigalg$-measurable due to continuity of $x(\cdot)$; the notion is therefore well-defined. An equivalent statement may be made up in terms of class-$\ClassKL$ functions: the system~\eqref{e:ssysdef} satisfies the strong global asymptotic stability in probability property (s-\gasp{}) if for every $\eta \in\;]0, 1[$ there exists a function $\beta\in\ClassKL$ such that $\Prob{\norm{x(t)}\le\beta(\norm{\xz}, t)\;\;\fa t\ge 0} \ge 1-\eta$. In the context of randomly switched systems this property can be derived from {\sc gas}~a.s.\ with the aid of the local Lipschitz property of the vector fields. We state this in the following proposition, whose proof is provided in~\secref{s:proofs:ss:gasp}.

		\begin{proposition}
			If~\eqref{e:ssysdef} is {\sc gas}~a.s., then it is s-\gasp{}.
			\label{p:gaspimplication}
		\end{proposition}

		In particular, the hypotheses of Theorem~\ref{t:gasasgen}, Theorem~\ref{t:gasaspoisson} and Theorem~\ref{t:gasasunif} each imply s-\gasp{} of~\eqref{e:ssysdef}.

	\section{Proofs of the Analysis Results}
	\label{s:proofs}
		The proofs of the theorems and corollaries of~\secref{s:mainres} and~\secref{s:general} are documented in this section. In order to simplify the presentation, a number of technical lemmas are stated and proved first in~\secref{s:lemmas}, followed by the proofs of the main results in~\secref{s:resproofs}. We carry out the proofs of Theorem~\ref{t:gasasunif} and Corollary~\ref{c:gasmunif}, both dealing with switching signals of class UH, in complete detail below. The proofs of Theorem~\ref{t:gasaspoisson} and Corollary~\ref{c:gasmpoisson} dealing with switching signals of class EH are similar and are sketched. We retain the notations and conventions of~\secref{s:prelims}. Let us recall some basic definitions and results.

		Let $I$ be a nonempty index set. A family of real-valued random variables $\{\xi_i\}_{i\in I}$ is said to be \emph{uniformly integrable}~\cite[Definition~3, p.~23]{ref:raoProbTheo} if
		\[
			\lim_{c\ra\infty}\sup_{i\in I}\Expec{\abs{\xi_i}\indic{\{\abs{\xi_i} > c\}}} = 0.
		\]
		The following Hadamard-de~la~Vall\'ee~Poussin criterion~\cite[Theorem~5, p.~24]{ref:raoProbTheo} for checking uniform integrability of a family of random variables will be employed later.

		\begin{proposition}[Hadamard-de~la~Vall\'ee~Poussin]
			A family of real-valued integrable random variables $\{\xi_i\}_{i\in I}$ is uniformly integrable if and only if there exists a convex function $\phi:\R\lra\posR$ with $\phi(0) = 0$ and $\lim_{r\ra\infty}\phi(r)/r = \infty$, such that $\sup_{i\in I}\Expec{\phi(\xi_i)} < \infty$.
			\label{p:poussin}
		\end{proposition}

		Recall that a family of random variables $(\xi_t)_{t\ge 0}$ converges almost surely (a.s.) if it converges pointwise outside a $\mathsf P$-null set. The following Proposition is standard, it can be readily derived from the Vitali convergence theorem~\cite[Theorem~4, p.~24]{ref:raoProbTheo}.

		\begin{proposition}
			\label{p:l1convfromas}
			If $(\xi_t)_{t\ge 0}$ is a \cadlag{} (i.e., right-continuous and possessing limits from the left) random process on the filtered probability space above, $(\xi_t)_{t\ge 0}$ is uniformly integrable, and $(\xi_t)_{t\ge 0}$ converges to $0$ a.s., then $\bigl(\Expec{\xi_t}\bigr)_{t\ge 0}$ converges to $0$.
		\end{proposition}

		We need Egorov's theorem on almost uniform convergence of a sequence of measurable functions (see e.g.,~\cite[Theorem~4, p.~50]{ref:raoProbTheo} for a proof).

		\begin{theorem}[Egorov]
			\label{t:egorov}
			Let $(g_n)_{n\in\N}$ be a sequence of measurable functions on $(\Omega, \sigalg, \mathsf P)$ and $g_n\ra g$ a.s. Then for every $\eps > 0$ there exists a measurable set $A_\eps$ with $\Prob{\Omega\setmin A_\eps} < \eps$ such that $\bigl(g_n\indic{A_\eps}\bigr)_{n\in\N}$ converges uniformly to $g\indic{A_\eps}$.
		\end{theorem}

		\subsection{Auxiliary lemmas}\mbox{}
		\label{s:lemmas}

		\begin{lemma}
			\label{l:ls}
			The system~\eqref{e:ssysdef} has the following property: for every $\eps > 0$ there exists $L_\eps > 0$ such that
			\begin{equation}
				\indic{]0, \eps[}(x(t))\left\lvert\frac{\drv\norm{x(t)}}{\drv t}\right\rvert \le L_\eps\norm{x(t)}.
				\label{e:asconv1}
			\end{equation}
			In particular, $\indic{]0, \eps[}(x(t))\norm{x(t)} \le \norm{\xz}\epower{L_\eps t}\;\;\fa t\ge 0$.
		\end{lemma}

		\begin{proof} Since $\{f_i\}_{i\in\PSet}$ is a finite family of locally Lipschitz vector fields, there exists some $\eps'' > 0$ and $L_{\eps''} > 0$ such that
			\[
				\sup_{\substack{i\in\PSet,\\\norm{x}\in[0, \eps''[}}\norm{f_i(x)} \le L_{\eps''}\norm{x}.
			\]
			Let $\eps := \eps'\wedge\eps''$. Note that $\fa x\in\R^n\setmin\{0\}$ we have 
			\[
				\left\lvert\frac{\drv\norm{x}^2}{\drv t}\right\rvert = \norm{2x\transp\frac{\drv x}{\drv t}} \le 2\norm{x}\norm{\frac{\drv x}{\drv t}}
			\]
			and
			\[
				\left\lvert\frac{\drv\norm{x}^2}{\drv t}\right\rvert = 2\norm{x}\abs{\frac{\drv\norm{x}}{\drv t}}.
			\]
			These two inequalities lead to $\left\lvert\frac{\drv\norm{x}}{\drv t}\right\rvert\le \norm{\frac{\drv x}{\drv t}}$. The inequality in~\eqref{e:asconv1} follows.
			Similarly, 
			\begin{equation}
				\frac{\drv\norm x}{\drv t} \le L_\eps\norm x\qquad \fa x\in\big\{x\in\R^n\,\big|\, \norm{x} < \eps\big\}\setmin\{0\}.
				\label{e:randstab:12}
			\end{equation}
			An application of a standard differential inequality~\cite[Theorem~1.2.1]{ref:lakshmikanthamDifferentialIntegralInequalities1} indicates that every solution $x(\cdot)$ of~\eqref{e:ssysdef} satisfies
			\[
				\norm{x(t)} \le \norm{\xz}\epower{L_\eps t}
			\]
			so long as $\norm{x(t)} < \eps$. This proves the claim.
		\end{proof}

		The following Barbalat-type lemma was stated without a complete proof in~\cite{ref:ransw}. It allows us to assert asymptotic convergence of $\norm{x(\cdot)}$ from the finiteness of a certain integral of $\norm{x(\cdot)}$.

		\begin{lemma}
			If $\alpha\in\ClassK$ and $\displaystyle{\int_0^\infty \alpha(\norm{x(t)})\drv t < \infty}$ a.s., then $\displaystyle{\lim_{t\ra\infty} \norm{x(t)} = 0}$ a.s., where $x(\cdot)$ is the solution of~\eqref{e:ssysdef}.
			\label{l:asconv}
		\end{lemma}

		\begin{proof}
			Suppose that the claim is false. Then there exists a measurable set $D$ of positive probability such that for every event in $D$ there exists some $\eps' > 0$ and a monotone increasing divergent sequence $(s_i)_{i\in\N}$ in $\posR$ such that $\alpha(\norm{x(s_i)}) > \eps'$ for all $i$.
			By the finiteness condition on the integral in the hypothesis, almost surely there exists $T(\eps) > 0$ such that 
			\begin{equation}
				\int_{T(\eps)}^\infty \alpha(\norm{x(t)})\drv t < \frac{1}{2}\int_0^{\frac{\ln 2}{L_\eps}}\alpha\biggl(\frac{\eps}{2}\epower{-L_\eps s}\biggr)\drv s,
				\label{e:asconv2}
			\end{equation}
			where the right hand side is a strictly positive quantity since $\alpha\in\ClassK$. For every event on a set of positive probability we have assumed that $(s_i)_{i\in\N}$ is a monotone increasing divergent sequence with $\alpha(\norm{x(s_i)}) > \eps$, and therefore there exists $i(\eps)\in\N$ such that $s_{i(\eps)} > T(\eps)$ with strictly positive probability. By continuity of $\norm{\cdot}$ and $x(\cdot)$, there exists an instant $t' > s_{i(\eps)}$ such that $\norm{x(t')} = \eps/2$, also with positive probability. But since $x(\cdot)$ solves~\eqref{e:ssysdef}, Lemma~\ref{l:ls} holds, and by~\eqref{e:asconv1} we have $\norm{x(t)}\in\;]0,\eps[$ for all $t\in\,]t', t'+\frac{\ln 2}{L_\eps}[$. Therefore
			\[
				\int_{t'}^{t'+\frac{\ln 2}{L_\eps}} \alpha(\norm{x(t)})\drv t \ge \int_{t'}^{t'+\frac{\ln 2}{L_\eps}}\alpha\biggl(\frac{\eps}{2}\epower{-L_\eps(t-t')}\biggr)\drv t
			\]
			with positive probability, which is a contradiction in view of~\eqref{e:asconv2}. The assertion follows.
		\end{proof}

		\begin{lemma}
			\label{l:Vatswtimesunif}
			Under the hypotheses of Theorem~\ref{t:gasasunif}, for each $j\in\N$ we have
			\[
			\Expec{V^{1+\kappa}_{\sigma(\tau_j)}(x(\tau_j))} \le \alpha_2^{1+\kappa}(\norm{\xz})\eta^j(\kappa),
			\]
			where $\displaystyle{\eta(\kappa) := \sum_{j\in\PSet}\frac{\mu^{1+\kappa}q_j\left(1-\epower{-\lambda_j(1+\kappa) T}\right)}{\lambda_j(1+\kappa)T}}$, $\kappa > 0$.
		\end{lemma}

		\begin{proof}
			Pick $i\in\Nz$. For $t\in[\tau_i, \tau_{i+1}[$, from (V2) we have
			\[
				V_{\sigma(\tau_{i+1})}(x(t)) \le V_{\sigma(\tau_{i+1})}(x(\tau_i))\epower{-\lambda_{\sigma(\tau_i)}(t-\tau_i)},
			\]
			and by continuity of $x(\cdot)$ and each Lyapunov function, and (V3),
			\[
				V_{\sigma(\tau_{i+1})}(x(t)) \le \mu V_{\sigma(\tau_i)}(x(\tau_i))\epower{-\lambda_{\sigma(\tau_i)}(t-\tau_i)}
			\]
			pointwise on $\Omega$. Fix $j\in\N$. For $\kappa > 0$, iterating the above inequality and employing the independence hypothesis (UH3) and (V1), we have
			\begin{equation}
			\begin{aligned}
				\Expec{V^{1+\kappa}_{\sigma(\tau_j)}(x(\tau_j))} & \le \alpha_2^{1+\kappa}(\norm{\xz})\Expec{\left(\prod_{i=0}^{j-1}\mu\epower{-\lambda_{\sigma(\tau_i)}S_{i+1}}\right)^{1+\kappa}}\\
				& = \alpha_2^{1+\kappa}(\norm{\xz})\prod_{i=0}^{j-1}\mu^{1+\kappa}\Expec{\epower{-\lambda_{\sigma(\tau_i)}(1+\kappa)S_{i+1}}}.
			\end{aligned}
			\label{e:Vatswtimesunif1}
			\end{equation}
			But
			\begin{align}
				\Expec{\epower{-\lambda_{\sigma(\tau_i)}(1+\kappa)S_{i+1}}} & = \Expec{\CExpec{\epower{-\lambda_{\sigma(\tau_i)}(1+\kappa)S_{i+1}}}{\sigalg_{\tau_i}}}\nonumber\\
				& = \Expec{\int_0^T \frac{1}{T} \epower{-\lambda_{\sigma(\tau_i)}(1+\kappa)s}\drv s}\nonumber\\
				& = \Expec{\frac{1-\epower{-\lambda_{\sigma(\tau_i)}(1+\kappa) T}}{\lambda_{\sigma(\tau_i)}(1+\kappa)T}}\nonumber\\
				& = \sum_{j\in\PSet}\frac{q_j\left(1-\epower{-\lambda_j(1+\kappa) T}\right)}{\lambda_j(1+\kappa)T}.
				\label{e:Vatswtimesunif2}
			\end{align}
			Substituting the right hand side of~\eqref{e:Vatswtimesunif2} in~\eqref{e:Vatswtimesunif1} leads to
			\[
			\Expec{V^{1+\kappa}_{\sigma(\tau_j)}(x(\tau_j))} \le \alpha_2^{1+\kappa}(\norm{\xz})\left(\sum_{i\in\PSet}\frac{\mu^{1+\kappa}q_i\left(1-\epower{-\lambda_i(1+\kappa) T}\right)}{\lambda_i(1+\kappa)T}\right)^{j},
			\]
			and considering the definition of $\eta(\kappa)$ the assertion follows.
		\end{proof}

		\begin{lemma}
			\label{l:finiteintunif}
			Under the hypotheses of Theorem~\ref{t:gasasunif} we have $\displaystyle{\int_0^\infty \!\alpha_1(\norm{x(t)})\drv t < \infty}$ a.s.
		\end{lemma}

		\begin{proof}
			For a fixed $t\in\posR$ we have
			\begin{align}
				\Expec{V_{\sigma(t)}(x(t))} & = \Expec{\sum_{i=0}^\infty V_{\sigma(t)}(x(t))\indic{\{t\in[\tau_i, \tau_{i+1}[\}}}\nonumber\\
				& = \sum_{i=0}^\infty \Expec{V_{\sigma(t)}(x(t))\indic{\{t\in[\tau_i, \tau_{i+1}[\}}},
				\label{e:finiteintunif1}
			\end{align}
			where we have employed the monotone convergence theorem~\cite[Theorem~1, \S1.3]{ref:raoProbTheo} to get the second equality. An application of (V1) and Tonelli's theorem~\cite[Theorem~11, \S1.3]{ref:raoProbTheo} gives us
			\begin{align}
				\Expec{\int_0^\infty \alpha_1(\norm{x(t)})\drv t} \le \Expec{\int_0^\infty V_{\sigma(t)}(x(t))\drv t} = \int_0^\infty \Expec{V_{\sigma(t)}(x(t))}\drv t,
				\label{e:finiteintunif2}
			\end{align}
			and in conjunction with~\eqref{e:finiteintunif1} we obtain
			\begin{align*}
				\Expec{\int_0^\infty \alpha_1(\norm{x(t)})\drv t} \le \int_0^\infty \sum_{i=0}^\infty \Expec{V_{\sigma(t)}(x(t))\indic{\{t\in[\tau_i, \tau_{i+1}[\}}}\drv t.
			\end{align*}
			A second application of monotone convergence theorem on the right hand side of the above leads to
			\begin{align*}
				\Expec{\int_0^\infty \alpha_1(\norm{x(t)})\drv t} & \le \sum_{i=0}^\infty \int_0^\infty \Expec{V_{\sigma(t)}(x(t))\indic{\{t\in[\tau_i, \tau_{i+1}[\}}}\drv t,
			\end{align*}
			and a further application of Tonelli's theorem on the right hand side gives
			\begin{align}
				 \sum_{i=0}^\infty \int_0^\infty \Expec{V_{\sigma(t)}(x(t))\indic{\{t\in[\tau_i, \tau_{i+1}[\}}}\drv t = \sum_{i=0}^\infty \Expec{\int_0^\infty V_{\sigma(t)}(x(t))\indic{\{t\in[\tau_i,\tau_{i+1}[\}}\drv t}.
				\label{e:finiteintunif3}
			\end{align}
			Each term in the series on the right hand side of~\eqref{e:finiteintunif3} may be estimated as follows:
			\begin{align*}
				\Expec{\int_0^\infty V_{\sigma(t)}(x(t))\indic{\{t\in[\tau_i, \tau_{i+1}[\}}\drv t} & \le \Expec{\int_0^\infty V_{\sigma(\tau_i)}(x(\tau_i))\epower{-\lambda_{\sigma(\tau_i)}(t-\tau_i)}\indic{\{t\in[\tau_i, \tau_{i+1}[\}}\drv t}
			\end{align*}
			by (V2), and therefore
			\begin{align}
				\Expec{\int_0^\infty V_{\sigma(t)}(x(t))\indic{\{t\in[\tau_i, \tau_{i+1}[\}}\drv t}
				& = \Expec{\int_{\tau_i}^{\tau_{i+1}} V_{\sigma(\tau_i)}(x(\tau_i))\epower{-\lambda_{\sigma(\tau_i)}(t-\tau_i)}\drv t}\nonumber\\
				& = \Expec{V_{\sigma(\tau_i)}(x(\tau_i))\left(\frac{1-\epower{-\lambda_{\sigma(\tau_i)}S_{i+1}}}{\lambda_{\sigma(\tau_i)}}\right)}\nonumber\\
				& = \Expec{\CExpec{V_{\sigma(\tau_i)}(x(\tau_i))\left(\frac{1-\epower{-\lambda_{\sigma(\tau_i)}S_{i+1}}}{\lambda_{\sigma(\tau_i)}}\right)}{\sigalg_{\tau_i}}}\nonumber\\
				& = \Expec{V_{\sigma(\tau_i)}(x(\tau_i))\left(\frac{1-\CExpec{\epower{-\lambda_{\sigma(\tau_i)}S_{i+1}}}{\sigalg_{\tau_i}}}{\lambda_{\sigma(\tau_i)}}\right)}\nonumber\\
				& = \Expec{\frac{V_{\sigma(\tau_i)}(x(\tau_i))}{\lambda_{\sigma(\tau_i)}}\left(1 - \int_0^T \frac{1}{T} \epower{-\lambda_{\sigma(\tau_i)}s}\drv s\right)}\nonumber\\
				& = \Expec{\frac{V_{\sigma(\tau_i)}(x(\tau_i))}{\lambda_{\sigma(\tau_i)}}\left(1-\frac{1-\epower{-\lambda_{\sigma(\tau_i)} T}}{\lambda_{\sigma(\tau_i)} T}\right)}\nonumber\\
				& \le M\Expec{V_{\sigma(\tau_i)}(x(\tau_i))},
				\label{e:finiteintunif4}
			\end{align}
			where $M := \max_{i\in\PSet}\left(\frac{1}{\lambda_i} - \frac{1-\epower{-\lambda_i T}}{\lambda_i^2 T}\right)$ is a well-defined positive real number because of the finiteness of $\PSet$. From~\eqref{e:finiteintunif3} and~\eqref{e:finiteintunif4} we get
			\begin{align*}
				\Expec{\int_0^\infty \alpha_1(\norm{x(t)})\drv t} & \le \sum_{i=0}^\infty \Expec{\int_0^\infty V_{\sigma(t)}(x(t))\indic{\{t\in[\tau_i,\tau_{i+1}[\}}\drv t}\\
				& \le M\alpha_2(\norm{\xz}) \sum_{i=0}^\infty \Expec{V_{\sigma(\tau_i)}(x(\tau_i))}\\
				& \le M\alpha_2(\norm{\xz}) \sum_{i=0}^\infty \eta^i(0)\\
				& < \infty,
			\end{align*}
			where $\eta$ is as defined in Lemma~\ref{l:Vatswtimesunif}, and $\eta(0)\in\;]0, 1[$ by (U3). This establishes the claim.
		\end{proof}

		\begin{lemma}
			Under the hypotheses of Theorem~\ref{t:gasasunif}, the family of random variables $\bigl\{V_{\sigma(t)}(x(t))\bigr\}_{t\ge 0}$ is uniformly integrable.
			\label{l:unifintegrgasmunif}
		\end{lemma}

		\begin{proof}
			To establish uniform integrability of the family $\bigl\{V_{\sigma(t)}(x(t))\bigr\}_{t\ge 0}$ we appeal to the Hadamard-de~la~Vall\'ee Poussin criterion in Proposition~\ref{p:poussin}. Since the function
			\[
			]-1, \infty[\;\ni r\mapsto \sum_{j\in\PSet}\frac{\mu^{1+r}q_j\left(1-\epower{-\lambda_j(1+r)T}\right)}{\lambda_j(1+r)T}\in\R
			\]
			is continuous, by (U3) there exists $\delta > 0$ such that $\sum_{j\in\PSet}\frac{\mu^{1+\delta} q_j\left(1-\epower{-\lambda_j(1+\delta) T}\right)}{\lambda_j(1+\delta)T} < 1$. The function $\phi(r) := r^{1+\delta}$ clearly is convex on $\posR$, and $\lim_{r\ra\infty}\phi(r)/r = \infty$. Let us prove that $\sup_{t\ge 0}\Expec{\bigl(V_{\sigma(t)}(x(t))\bigr)^{1+\delta}} <  \infty$.

			First let us note that for each $i\in\Nz$ the function $V^{1+\delta}_{\sigma(t)}(x(t))\indic{\{t\in[\tau_i, \tau_{i+1}[\}}$ is integrable for arbitrary $t\in\posR$. Indeed,
			\begin{align*}
				\Expec{V^{1+\delta}_{\sigma(t)}(x(t))\indic{\{t\in[\tau_i, \tau_{i+1}[\}}} & = \Expec{\CExpec{V^{1+\delta}_{\sigma(t)}(x(t))\indic{\{t\in[\tau_i, \tau_{i+1}[\}}}{\sigalg_{\tau_i}}}\nonumber\\
				& \le \Expec{\CExpec{V^{1+\delta}_{\sigma(\tau_i)}(x(\tau_i))\epower{-\lambda_{\sigma(\tau_i)}(1+\delta)(t-\tau_i)}\indic{\{t\in[\tau_i, \tau_{i+1}[\}}}{\sigalg_{\tau_i}}}\nonumber\\
				& = \Expec{V^{1+\delta}_{\sigma(\tau_i)}(x(\tau_i))\epower{-\lambda_{\sigma(\tau_i)}(1+\delta)(t-\tau_i)}\CExpec{\indic{\{t\in[\tau_i, \tau_{i+1}[\}}}{\sigalg_{\tau_i}}},
			\end{align*}
			and since $S_{i+1}$ is uniform-$T$ and independent of $\sigalg_{\tau_i}$, we have 
			\begin{align*}
				\CExpec{\indic{\{t\in[\tau_i, \tau_{i+1}[\}}}{\sigalg_{\tau_i}} & = \indic{\{t\in[\tau_i,\infty[\}}\CProb{S_{i+1} > t-\tau_i}{\sigalg_{\tau_i}}\\
				& = \left(\left(1-\frac{t-\tau_i}{T}\right)\mx 0\right).
			\end{align*}
			Therefore,
			\begin{align}
				& \Expec{V^{1+\delta}_{\sigma(t)}(x(t))\indic{\{t\in[\tau_i, \tau_{i+1}[\}}}\nonumber\\
				& \quad \le \Expec{V^{1+\delta}_{\sigma(\tau_i)}(x(\tau_i))\epower{-\lambda_{\sigma(\tau_i)}(1+\delta)(t-\tau_i)}\left(\left(1-\frac{t-\tau_i}{T}\right)\mx 0\right)\indic{\{t\in[\tau_i, \infty[\}}}.
				\label{e:gasmunif1}
			\end{align}
			By definition of $\delta$, the right hand side of~\eqref{e:gasmunif1} is at most $M\Expec{V^{1+\delta}_{\sigma(\tau_i)}(x(\tau_i))}$, where $M := \exp\bigl(\min_{j\in\PSet}\lambda_j\cdot (1+\delta)T\bigr)$. Lemma~\ref{l:Vatswtimesunif} with $\kappa = \delta$ shows that 
			\begin{equation}
				\Expec{V^{1+\delta}_{\sigma(\tau_i)}(x(\tau_i))} \le \alpha_2^{1+\delta}(\norm{\xz})\eta(\delta)^i,
				\label{e:gasmunif2}
			\end{equation}
			where $\eta(\delta) \in\;]0, 1[$ by construction. By~\eqref{e:gasmunif1} we know that the random variable $V^{1+\delta}_{\sigma(t)}(x(t))\indic{\{t\in[\tau_i, \tau_{i+1}[\}}$ is integrable for each $i$; we can therefore apply the monotone convergence theorem to arrive at
			\begin{align}
				\Expec{\bigl(V_{\sigma(t)}(x(t))\bigr)^{1+\delta}} & = \Expec{\left(\sum_{i=0}^\infty V_{\sigma(t)}(x(t))\indic{\{t\in[\tau_i, \tau_{i+1}[\}}\right)^{1+\delta}}\nonumber\\
				& = \Expec{\sum_{i=0}^\infty V^{1+\delta}_{\sigma(t)}(x(t))\indic{\{t\in[\tau_i, \tau_{i+1}[\}}}\nonumber\\
				& = \sum_{i=0}^\infty \Expec{V^{1+\delta}_{\sigma(t)}(x(t))\indic{\{t\in[\tau_i, \tau_{i+1}[\}}}.
				\label{e:gasmunif3}
			\end{align}
			We know from~\eqref{e:gasmunif2} that $\Expec{V^{1+\delta}_{\sigma(t)}(x(t))\indic{\{t\in[\tau_i, \tau_{i+1}[\}}} \le M\alpha_2^{1+\delta}(\norm{\xz})\eta^i(\delta)$ for each $i\in\Nz$. Substitution in~\eqref{e:gasmunif3} leads to
			\begin{equation}
			\begin{aligned}
				\sup_{t\ge 0}\Expec{\bigl(V_{\sigma(t)}(x(t))\bigr)^{1+\delta}} & = \sup_{t\ge 0}\sum_{i=0}^\infty \Expec{V^{1+\delta}_{\sigma(t)}(x(t))\indic{\{t\in[\tau_i, \tau_{i+1}[\}}}\\
				& \le \sup_{t\ge 0} M\alpha_2^{1+\delta}(\norm{\xz})\sum_{i=0}^\infty \eta^i(\delta)\\
				& < \infty.
			\end{aligned}
			\label{e:supVsigmatunif}
			\end{equation}
			This shows that the family $\bigl\{V_{\sigma(t)}(x(t))\bigr\}_{t\ge 0}$ is uniformly integrable.
		\end{proof}

		\begin{lemma}
			\label{l:geniterate}
			Under the hypotheses of Theorem~\ref{t:gasasgen}, for every $\nu\in\N$ we have $\Expec{V_{\sigma(\tau_{\nu})}(x(\tau_{\nu}))} \le \theta^{\nu} V_{\sigmaz}(\xz)$.
		\end{lemma}

		\begin{proof}
			Fix $i\in\Nz$. For $t\in[\tau_i, \tau_{i+1}[$ and $j\in\PSet$, from (V2$'$) we have
			\[
				V_j(x(t)) \le V_j(x(\tau_i))\epower{-\lambda_{j\sigma(\tau_i)}(t-\tau_i)}.
			\]
			In particular, for $t\in[\tau_i, \tau_{i+1}[$,
			\[
				V_{\sigma(\tau_{i+1})}(x(t)) \le V_{\sigma(\tau_{i+1})}(x(\tau_i))\epower{-\lambda_{\sigma(\tau_{i+1}),\sigma(\tau_i)}(t-\tau_i)},
			\]
			and by continuity of $x(\cdot)$ and each Lyapunov function, and (V3$'$),
			\[
				V_{\sigma(\tau_{i+1})}(x(t)) \le \mu V_{\sigma(\tau_i)}(x(\tau_i))\epower{-\lambda_{\sigma(\tau_{i+1}), \sigma(\tau_i)}(t-\tau_i)}
			\]
			pointwise on $\Omega$. Therefore,
			\begin{equation}
				\label{e:condViterategen}
				\CExpec{V_{\sigma(\tau_{i+1})}(x(\tau_{i+1}))}{\sigalg_{\tau_i}} \le \mu V_{\sigma(\tau_i)}(x(\tau_i))\CExpec{\epower{-\lambda_{\sigma(\tau_{i+1}), \sigma(\tau_i)}S_{i+1}}}{\sigalg_{\tau_i}}.
			\end{equation}
			(GH3) shows that $S_{i+1}$ and $\sigma(\tau_{i+1})$ are conditionally independent given $\sigalg_{\tau_i}$, and therefore,
			\[
				\CExpec{\epower{-\lambda_{\sigma(\tau_{i+1}),\sigma(\tau_i)}S_{i+1}}}{\sigalg_{\tau_i}} = \sum_{j\in\PSet}\CExpec{\epower{-\lambda_{j,\sigma(\tau_i)}S_{i+1}}}{\sigalg_{\tau_i}} p_{\sigma(\tau_i), j}.
			\]
			Since $\sigma(\tau_i)$ is $\sigalg_{\tau_i}$-measurable,
			\[
				\sum_{j\in\PSet}\CExpec{\epower{-\lambda_{j,\sigma(\tau_i)}S_{i+1}}}{\sigalg_{\tau_i}} p_{\sigma(\tau_i), j} \le \max_{k\in\PSet}\sum_{j\in\PSet}\Expec{\epower{-\lambda_{j,k}S_1}} p_{k, j}.
			\]
			By (G3) there exists a $\theta\in\:]0, 1[$ such that the quantity on the right hand side of the above inequality is at most $\theta/\mu$. Therefore, we get
			\[
				\mu \CExpec{\epower{\lambda_{\sigma(\tau_{i+1}), \sigma(\tau_i)}S_{i+1}}}{\sigalg_{\tau_i}} \le \theta < 1,
			\]
			which in view of~\eqref{e:condViterategen} shows that
			\[
				\CExpec{V_{\sigma(\tau_{i+1})}(x(\tau_{i+1}))}{\sigalg_{\tau_i}} \le \theta V_{\sigma(\tau_i)}(x(\tau_i)).
			\]
			Fixing $\nu\in\N$, since $(\tau_i)_{i\in\N}$ is an increasing sequence of $(\sigalg_t)_{t\ge 0}$-optional times, it follows from standard properties of conditional expectations\footnote{The property being utilized is the following: If $\tau$ and $\tau'$ are $(\sigalg_t)_{t\ge 0}$-optional times, and $\tau \le \tau'$, then $\sigalg_{\tau}$ is a sub-sigma-algebra of $\sigalg_{\tau'}$. See e.g.,~\cite[Chapter~6]{ref:raoProbTheo} for further details.}
			\begin{align*}
				\Expec{V_{\sigma(\tau_\nu)}(x(\tau_\nu))} & = \Expec{\CExpec{\cdots\CExpec{\CExpec{V_{\sigma(\tau_\nu)}(x(\tau_\nu))}{\sigalg_{\tau_{\nu-1}}}}{\sigalg_{\tau_{\nu-2}}}\cdots}{\sigalg_{\tau_1}}}\\
				& \le \Expec{\CExpec{\cdots\CExpec{\theta V_{\sigma(\tau_{\nu-1})}(x(\tau_{\nu-1})) }{\sigalg_{\tau_{\nu-2}}}\cdots}{\sigalg_{\tau_1}}}\\
				& \le \theta^\nu V_{\sigmaz}(\xz).
			\end{align*}
			This proves the assertion.
		\end{proof}

		\begin{lemma}
			\label{l:finiteintgen}
			Under the hypotheses of Theorem~\ref{t:gasasgen} we have $\displaystyle{\int_0^\infty \alpha_1(\norm{x(t)}) \drv t < \infty}$ a.s.
		\end{lemma}

		\begin{proof}
			Following the proof of Lemma~\ref{l:finiteintunif} we have
			\[
				\Expec{V_{\sigma(t)}(x(t))} = \sum_{i=0}^\infty \Expec{V_{\sigma(t)}(x(t))\indic{\{t\in[\tau_i, \tau_{i+1}[\}}}.
			\]
			From (V1$'$), the monotone convergence theorem, and two applications of Tonelli's theorem, (as in the proof of Lemma~\ref{l:finiteintunif},) we get
			\begin{align}
				\label{e:gensum}
				\Expec{\int_0^\infty \alpha_1(\norm{x(t)})\drv t} \le \int_0^\infty\Expec{V_{\sigma(t)}(x(t))\drv t} = \sum_{i=0}^\infty \Expec{\int_{\tau_i}^{\tau_{i+1}} V_{\sigma(t)}(x(t))\drv t}.
			\end{align}
			Now by (V2$'$) we get
			\begin{align*}
				\Expec{\int_{\tau_i}^{\tau_{i+1}} V_{\sigma(t)}(x(t)) \drv t} & \le \Expec{V_{\sigma(\tau_i)}(x(\tau_i))\CExpec{\int_{\tau_i}^{\tau_{i+1}} \epower{-\lambda_{\sigma(\tau_i),\sigma(\tau_i)}(t-\tau_i)}\drv t}{\sigalg_{\tau_i}}}\\
				& = \Expec{V_{\sigma(\tau_i)}(x(\tau_i)) \left(\frac{1-\CExpec{\epower{-\lambda_{\sigma(\tau_i),\sigma(\tau_i)}S_{i+1}}}{\sigalg_{\tau_i}}}{\lambda_{\sigma(\tau_i),\sigma(\tau_i)}}\right)}.
			\end{align*}
			Note that the non-degeneracy of the matrix $Q$ yields $\Expec{\epower{-\lambda_{i, i}S_1}} < \infty$ for all $i\in\PSet$. This together with the fact that $\sigma(\tau_i)$ is $\sigalg_{\tau_i}$-measurable, guarantees the existence of a constant $M > 0$, such that
			\[
				\Expec{\int_{\tau_i}^{\tau_{i+1}} V_{\sigma(t)}(x(t))\drv t} \le M\Expec{V_{\sigma(\tau_i)}(x(\tau_i))}.
			\]
			Substituting in~\eqref{e:gensum} we arrive at
			\[
				\Expec{\int_0^\infty \alpha_1(\norm{x(t)})\drv t} \le \sum_{i=0}^\infty M \Expec{V_{\sigma(\tau_i)}(x(\tau_i))} \le M \alpha_2(\norm{\xz})\sum_{i=0}^\infty \theta^i < \infty
			\]
			in view of Lemma~\ref{l:geniterate} and (V3$'$). We immediately get $\Prob{\int_0^\infty \alpha_1(\norm{x(t)})\drv t < \infty} = 1$, as asserted.
		\end{proof}

		\subsection{Proofs of the Results in~\secref{s:mainres} and~\secref{s:general}}\mbox{}
		\label{s:resproofs}
		As stated at the beginning of~\secref{s:proofs}, the proofs of Theorem~\ref{t:gasasunif} and Corollary~\ref{c:gasmunif} are carried out in detail below, following which we provide sketches of the proofs of Theorem~\ref{t:gasaspoisson} and Corollary~\ref{c:gasmpoisson}.

		\emph{Proof of Theorem~\ref{t:gasasunif}.}
			To see the property (AS2) of~\eqref{e:ssysdef} we note that by Lemma~\ref{l:finiteintunif}, $\Prob{\int_0^\infty \alpha_1(\norm{x(t)})\drv t < \infty} = 1$. Lemma~\ref{l:asconv} now shows that $\norm{x(t)} \ra 0$ a.s.\ as $t\ra\infty$ since $\alpha_1\in\ClassKinfty$. Since $\xz$ was arbitrary, to establish (AS2) it only remains to show that the solutions corresponding to all initial conditions $\xz'$ with $\norm{\xz'} < \norm{\xz}$ are also asymptotically convergent. To this end, observe that for every fixed $\omega\in\Omega$, $\nu\in\N$, and $t\in[\tau_\nu(\omega), \tau_{\nu+1}(\omega)[$, a straightforward computation with the aid of (V1)-(V3) gives
			\begin{equation}
				\label{e:comparunif}
				V_{\sigma(t, \omega)}(x(t, \omega)) \le \alpha_2(\norm{\xz})\mu^\nu \prod_{i=0}^{\nu-1} \epower{-\lambda_{\sigma(\tau_{i}(\omega), \omega)}S_{i+1}(\omega)}\epower{-\lambda_{\sigma(\tau_\nu(\omega), \omega)} (t-\tau_\nu(\omega))}.
			\end{equation}
			Here $x(\cdot, \omega)$ corresponds to the solution of~\eqref{e:ssysdef} initialized at $\xz$. If $x'(\cdot, \omega)$ denotes the solution corresponding to the initial condition $\xz'$, then from~\eqref{e:comparunif} we have
			\[
				V_{\sigma(t, \omega)}(x'(t, \omega)) < \alpha_2(\norm{\xz})\mu^\nu \prod_{i=0}^{\nu-1} \epower{-\lambda_{\sigma(\tau_{i}(\omega), \omega)}S_{i+1}(\omega)}\epower{-\lambda_{\sigma(\tau_\nu(\omega), \omega)} (t-\tau_\nu(\omega))}
			\]
			whenever $\norm{\xz'} < \norm{\xz}$, since the right-hand side of~\eqref{e:comparunif} depends on the initial condition only through the function $\alpha_2$, which is monotone increasing. This proves (AS2).

			Now we verify (AS1). Fix $\eps > 0$. We know from the (AS2) property proved above that almost surely there exists $T(1, \eps) > 0$ such that $\norm{\xz} < 1$ implies that $\sup_{t\ge T(1,\eps)}\norm{x(t)} < \eps$. Select $\delta(\eps) = \min\left\{\eps \epower{-L_\eps T(1, \eps)}, 1\right\}$. By Lemma~\ref{l:ls}, $\norm{\xz} < \delta(\eps)$ implies
			\[
				\norm{x(t)} \le \norm{\xz}\epower{L_\eps t} < \delta(\eps)\epower{L_\eps T(1, \eps)} < \eps \quad \fa t\in[0, T(1, \eps)].
			\]
			Further, the (AS2) property guarantees that with the above choice of $\delta$ and $\xz$, we have $\sup_{t\ge T(1, \eps)} \norm{x(t)} < \eps$ for events in a set of full measure. Thus, $\norm{\xz} < \delta(\eps)$ implies that $\sup_{t \ge 0}\norm{x(t)} < \eps$ a.s. Since $\eps$ is arbitrary, the (AS1) property of~\eqref{e:ssysdef} follows.

			We conclude that~\eqref{e:ssysdef} is \gasas{}$\;\;$\endproof

		\emph{Proof of Theorem~\ref{t:gasaspoisson} (Sketch).}
			First we observe that under the hypotheses of Theorem~\ref{t:gasaspoisson}, for each $j\in\N$ we have
			\[
			\Expec{V^{1+\kappa}_{\sigma(\tau_j)}(x(\tau_j))} \le \alpha_2^{1+\kappa}(\norm{\xz})\eta^j(\kappa)\quad\text{whenever $(1+\kappa)\lambda_i + \lambda > 0$ for all $i\in\PSet$},
			\]
			where $\displaystyle{\eta(\kappa) := \sum_{j\in\PSet}\frac{\mu^{1+\kappa} q_j}{1+\lambda_j(1+\kappa)/\lambda}}$, $\kappa > 0$. This can be proved along the lines of Lemma~\ref{l:Vatswtimesunif}. In particular, at the step corresponding to~\eqref{e:Vatswtimesunif2} we employ the (E3) condition $(1+\kappa)\min_{i\in\PSet}\lambda_i + \lambda > 0$ as
			\begin{align*}
				\Expec{\epower{-\lambda_{\sigma(\tau_i)}(1+\kappa)S_{i+1}}} & = \Expec{\CExpec{\epower{-\lambda_{\sigma(\tau_i)}(1+\kappa)S_{i+1}}}{\sigalg_{\tau_i}}}\\
				& = \Expec{\lambda \int_0^\infty \epower{-\bigl(\lambda_{\sigma(\tau_i)}(1+\kappa) + \lambda\bigr)s}\drv s}\\
				& = \sum_{j\in\PSet}\frac{q_j}{1+(1+\kappa)\lambda_j/\lambda}.
			\end{align*}
			Second we observe that $\displaystyle{\int_0^\infty \!\alpha_1(\norm{x(t)})\drv t < \infty}$ a.s. The proof is similar to that of Lemma~\ref{l:finiteintunif}; the only difference lies in the step corresponding to~\eqref{e:finiteintunif4}, where we employ the condition $(1+\kappa)\min_{i\in\PSet}\lambda_i + \lambda > 0$ to arrive at
			\[
				\Expec{\int_0^\infty V_{\sigma(t)}(x(t))\indic{\{t\in[\tau_i, \tau_{i+1}[\}}\drv t} \le \Expec{V_{\sigma(\tau_i)}(x(\tau_i))}\frac{1}{\min_{j\in\PSet}\lambda_j + \lambda}.
			\]
			The subsequent steps follow those of Lemma~\ref{l:finiteintunif} and we get
			\[
				\Expec{\int_0^\infty \alpha_1(\norm{x(t)})\drv t} \le \frac{\alpha_2(\norm{\xz})}{\min_{j\in\PSet}\lambda_j + \lambda} \sum_{i=0}^\infty \eta^j(0) < \infty,
			\]
			where $\eta$ is as defined at the beginning of the current proof. With these ingredients, to see the property (AS2) of~\eqref{e:ssysdef} we note that in view of $\Prob{\int_0^\infty \alpha_1(\norm{x(t)})\drv t < \infty} = 1$, Lemma~\ref{l:asconv} gives $\norm{x(t)} \ra 0$ a.s.\ as $t\ra\infty$ since $\alpha_1\in\ClassKinfty$. This proves (AS2) because the only dependence on the initial condition is through $\alpha_2(\norm{\xz})$ and $\xz$ is arbitrary (as argued in the Proof of Theorem~\ref{t:gasasunif} above). The proof of (AS1) is identical to that in the proof of Theorem~\ref{t:gasasunif}, and we omit the details. It follows that~\eqref{e:ssysdef} is \gasas{}$\;\;$\endproof

		\emph{Proof of Corollary~\ref{c:gasmunif}.}
			Our first objective is to prove asymptotic convergence of the net $\bigl(\Expec{\alpha_1(\norm{x(t)})}\bigr)_{t\ge 0}$ to $0$. We have proved global asymptotic convergence a.s.\ of the process $(x(t))_{t\ge 0}$ to $0$ in Theorem~\ref{t:gasasunif}, and via hypothesis (V1) this shows that the process $\bigl(V_{\sigma(t)}(x(t))\bigr)_{t\ge 0}$ also converges a.s.\ to $0$ since $\alpha_2\in\ClassKinfty$. From Lemma~\ref{l:unifintegrgasmunif} we know that the family $\bigl\{V_{\sigma(t)}(x(t))\bigr\}_{t\ge 0}$ is uniformly integrable, and by Proposition~\ref{p:l1convfromas} it follows that $\lim_{t\ra\infty}\Expec{V_{\sigma(t)}(x(t))} = 0$. This implies global asymptotic convergence of $\Expec{\alpha_1(\norm{x(t)})}$ to $0$ in the light of (V1), and verifies the (SM2) property with $\alpha = \alpha_1$.

			It remains to prove (SM1). Following the notation of the proof of Lemma~\ref{l:unifintegrgasmunif}, we note that $\eta(0)\in\;]0, 1[$ by (U3). To establish (SM1) we only need to note that with $\delta = 0$ in~\eqref{e:supVsigmatunif} we have
			\[
				\sup_{t\ge 0}\Expec{V_{\sigma(t)}(x(t))} \le M\alpha_2(\norm{\xz})\frac{1}{1-\eta(0)}.
			\]
			For $\eps > 0$ preassigned, we choose $\wt\delta < \alpha_2^{-1}\bigl(\eps(1+\eta(0))/M\bigr)$ to see that 
			\[
				\sup_{t\ge 0}\Expec{\alpha_1(\norm{x(t)})} < \eps \quad\text{whenever }\norm{\xz} < \wt\delta.
			\]
			The (SM1) property with $\alpha = \alpha_1$ follows, thereby completing the proof.$\;\;$\endproof

		\emph{Proof of Corollary~\ref{c:gasmpoisson} (Sketch).}
			We follow the proof of Corollary~\ref{c:gasmunif} above. Since the proof of (SM1) is identical to that in the aforesaid proof, we give the details for the proof of (SM2). This involves establishing asymptotic convergence of the net $\bigl(\Expec{\alpha_1(\norm{x(t)})}\bigr)_{t\ge 0}$ to $0$. Since global asymptotic convergence of the process $(x(t))_{t\ge 0}$ to $0$ has been established in Theorem~\ref{t:gasaspoisson}, in the light of (V1) and Proposition~\ref{p:l1convfromas} it suffices to show that the family $\bigl\{V_{\sigma(t)}(x(t))\bigr\}_{t\ge 0}$ is uniformly integrable to conclude that $\lim_{t\ra\infty} \Expec{V_{\sigma(t)}(x(t))} = 0$.

			To this end, we need to follow the steps of Lemma~\ref{l:unifintegrgasmunif} above to establish uniform integrability of $\bigl\{V_{\sigma(t)}(x(t))\bigr\}_{t\ge 0}$. Since the function $]-1, \infty[\;\ni r\mapsto (1+r)\lambda_i + \lambda\in\R$ is continuous for each $i\in\PSet$ and $\PSet$ is a finite set, by (E3) there exists $\delta' > 0$ such that $(1+\delta')\lambda_i + \lambda > 0$ for all $i\in\PSet$. Also, since the function
			\[
				]-1, \infty[\;\ni r\mapsto \sum_{j\in\PSet}\frac{\mu^{1+r} q_j}{1+(1+r)\lambda_j/\lambda}\in\R
			\]
			is continuous, by (E4) there exists $\delta'' > 0$ such that $\sum_{j\in\PSet}\frac{\mu^{1+\delta''} q_j}{1+(1+\delta'')\lambda_j/\lambda} < 1$. Let $\delta := \delta'\mn\delta''$. The function $\phi(r) := r^{1+\delta}$ clearly is convex on $\posR$, and $\lim_{r\ra\infty}\phi(r)/r = \infty$. If we prove that $\sup_{t\ge 0}\Expec{\bigl(V_{\sigma(t)}(x(t))\bigr)^{1+\delta}} <  \infty$, then the Hadamard-de~la~Vall\'ee Poussin criterion in Proposition~\ref{p:poussin} may be applied to conclude uniform integrability of $\bigl\{V_{\sigma(t)}(x(t))\bigr\}_{t\ge 0}$.

			Calculations show that the inequality corresponding to~\eqref{e:gasmunif1} can be written as
			\begin{align*}
				\Expec{V^{1+\delta}_{\sigma(t)}(x(t))\indic{\{t\in[\tau_i, \tau_{i+1}[\}}} & \le \Expec{V^{1+\delta}_{\sigma(\tau_i)}(x(\tau_i))\epower{-(\lambda_{\sigma(\tau_i)}(1+\delta) + \lambda)(t-\tau_i)}\indic{\{t\in[\tau_i, \infty[\}}},
			\end{align*}
			and that corresponding to~\eqref{e:gasmunif2} can be written as
			\[
				\Expec{V^{1+\delta}_{\sigma(\tau_i)}(x(\tau_i))} \le \alpha_2^{1+\delta}(\norm{\xz})\eta(\delta)^i,
			\]
			where $\displaystyle{\eta(\kappa) := \sum_{j\in\PSet}\frac{\mu^{1+\kappa} q_j}{1+\lambda_j(1+\kappa)/\lambda}}$. The step corresponding to~\eqref{e:gasmunif3} is identical, and the one corresponding to~\eqref{e:supVsigmatunif} is
			\begin{align*}
				\sup_{t\ge 0}\Expec{\bigl(V_{\sigma(t)}(x(t))\bigr)^{1+\delta}} & = \sup_{t\ge 0}\sum_{i=0}^\infty \Expec{V^{1+\delta}_{\sigma(t)}(x(t))\indic{\{t\in[\tau_i, \tau_{i+1}[\}}}\\
				& \le \sup_{t\ge 0} \alpha_2^{1+\delta}(\norm{\xz})\sum_{i=0}^\infty \eta^i(\delta)\\
				& < \infty.
			\end{align*}
			This concludes the proof.$\;\;$\endproof

		\emph{Proof of Theorem~\ref{t:gasasgen}.}
			The proof mimics that of Theorem~\ref{t:gasasunif} above; the only change required here is to replace the occurrence of Lemma~\ref{l:finiteintunif} by Lemma~\ref{l:finiteintgen}.$\;\;$\endproof

		\subsection{Proof of Proposition~\ref{p:gaspimplication}}\mbox{}
		\label{s:proofs:ss:gasp}
		\emph{Proof of Proposition~\ref{p:gaspimplication}.}
			Let us verify property (ii) of Definition~\ref{d:gasp} assuming that~\eqref{e:ssysdef} is \gasas{} Fix $\eta, r, \eps' > 0$ and $\xz\in\R^n$ with $\norm{\xz} < r$. Since $\{f_i\}_{i\in\PSet}$ is a finite set of locally Lipschitz vector fields, there exists $L_{\eps'} > 0$ such that $\sup_{i\in\PSet,\norm{x} < \eps'}\norm{f_i(x)} \le L_{\eps'}\norm x$. Let $c := \frac{\ln 2}{L_{\eps'}}$, and define the sequence of time instants $(s_j)_{j\in\Nz}$ such that $s_0 := 0$ and $s_j - s_{j-1} = c$ for every $j\in\N$. By the (AS2) property of~\eqref{e:ssysdef} we have $\Prob{\lim_{t\ra\infty}\norm{x(t)} = 0} = 1$, which also implies that $\Prob{\lim_{i\ra\infty} \norm{x(s_i)} = 0} = 1$. By Egorov's Theorem~\ref{t:egorov} there exists a measurable set $A_{\eta}$ such that $\Prob{\Omega\setmin A_{\eta}} < \eta$ and $\bigl(x(s_i)\indic{A_\eta}\bigr)_{i\in\N}$ uniformly converges to $0$. The uniform convergence condition by definition implies that there exists $i_0\in\N$ such that $\sup_{i\ge i_0}\bigl(\norm{x(s_i)}\indic{A_{\eta}}\bigr) < \frac{\eps'}{2}$. By construction of the sequence $(s_i)_{i\in\N}$ we must have $\norm{x(t)}\indic{A_{\eta}} < \eps'$ for all $t\ge s_{i_0}$ in view of continuity of $x(\cdot)$. To see this, fix a time $t' > s_{i_0}$. The construction of the sequence $(s_i)_{i\in\N}$ shows that there exists a $j(t')\in\N$ such that $t'\in[s_{j(t') - 1}, s_{j(t')}[$. The local Lipschitz condition on the set of vector fields $\{f_i\}_{i\in\PSet}$ implies that
			\[
			\norm{x(t')}\indic{A_\eta} \le \sup_{s\in[s_{j(t')-1}, s_{j(t')}[}\norm{x(s)}\indic{A_\eta} < \frac{\eps'}{2}\epower{L_{\eps'}\left(s-s_{j(t')}\right)} < \frac{\eps'}{2}\epower{L_{\eps'}c} = \eps',
			\]
			where the last equality is true by definition of $c$. Since $t'$ was arbitrary, the assertion follows. Since $\xz$ was arbitrary, to establish the property (ii) of Definition~\ref{d:gasp} it only remains to show that the solutions restricted to $A_\eta$ corresponding to all initial conditions $\xz'$ with $\norm{\xz'} < \norm{\xz}$ are also asymptotically convergent. To this end, observe that for every fixed $\omega\in\Omega$, and therefore for every fixed $\omega\in A_\eta$, $\nu\in\N$, and $t\in[\tau_\nu(\omega), \tau_{\nu+1}(\omega)[$, a straightforward computation with the aid of (V1)-(V3) gives
			\begin{equation}
				\label{e:compar}
				V_{\sigma(t, \omega)}(x(t, \omega)) \le \alpha_2(\norm{\xz})\mu^\nu \prod_{i=0}^{\nu-1} \epower{-\lambda_{\sigma(\tau_{i}(\omega), \omega)}S_{i+1}(\omega)}\epower{-\lambda_{\sigma(\tau_\nu(\omega), \omega)}(t-\tau_\nu(\omega))}.
			\end{equation}
			Here $x(\cdot, \omega)$ corresponds to the solution of~\eqref{e:ssysdef} initialized at $\xz$. If $x'(\cdot, \omega)$ denotes the solution corresponding to the initial condition $\xz'$, then from~\eqref{e:compar} we have
			\[
				V_{\sigma(t, \omega)}(x'(t, \omega)) < \alpha_2(\norm{\xz})\mu^\nu \prod_{i=0}^{\nu-1} \epower{-\lambda_{\sigma(\tau_{i}(\omega), \omega)}S_{i+1}(\omega)}\epower{-\lambda_{\sigma(\tau_\nu(\omega), \omega)}(t-\tau_\nu(\omega))}
			\]
			whenever $\norm{\xz'} < \norm{\xz}$, since the right-hand side of~\eqref{e:compar} depends on the initial condition only through the function $\alpha_2$, which is monotone increasing. This proves (ii). To establish (i), let us fix $\eta\in\;]0, 1[$ and $\eps > 0$. By (ii) there exists a $T > 0$ corresponding to $\eta' = \eta$, $r = 1$ and $\eps' = \eta$ such that $\norm{\xz} < 1$ implies that $\sup_{t\ge T}\norm{x(t)}\indic{A_\eta} < \eps$. The local Lipschitz condition on the set of vector fields $\{f_i\}_{i\in\N}$ guarantees the existence of a positive $\delta' > 0$ such that $\sup_{t\in[0, T]}\norm{x(t)} < \eps$ whenever $\norm{\xz} < \delta$. Picking $\delta = 1\mn\delta'$ we see that $\norm{\xz} < \delta$ implies that $\sup_{t\ge 0}\norm{x(t)}\indic{A_\eta} < \eps$. The implication is now completely established.$\;\;$\endproof

	\section{Control Synthesis}
		\label{s:syn}
		Our goal in this section is to synthesize feedback control functions for stabilization (in a suitable stochastic sense) of randomly switched systems with control inputs. For brevity, we shall restrict ourselves to controllers which render the closed-loop switched system {\sc gas}~a.s.\ for a switching signal of class EH. The results automatically give the $\alpha_1$-\gasm{} property also in addition to {\sc gas}~a.s., in view of the close relationship between the sufficient conditions for {\sc gas}~a.s.\ and \gasm{} in our analysis results of~\secref{s:mainres}.

		There are two distinct and obvious controller architectures: one in which the the control function depends on the switching signal $\sigma$, and the other in which the control function does not depend on $\sigma$. In the first case, which is presented in~\secref{s:synmd}, we combine universal formulae for feedback stabilization of nonlinear systems with our analysis results to design controllers which ensure {\sc gas}~a.s.\ of the closed-loop switched system. In the second case, which is presented in~\secref{s:synmi}, we search for a controller which stabilizes some subsystems while not destabilizing the others too much, and with the aid of our analysis results, ensure that the closed-loop switched system is \gasas{}

		\subsection{Mode-dependent Controllers}
		\label{s:synmd}
		Consider the affine-in-control switched system
		\begin{equation}
			\label{e:ssysdefcon}
			\dot x = f_\sigma(x) + \sum_{j = 1}^mg_{\sigma, j}(x) u_j,\qquad (x(0), \sigma(0)) = (\xz,\sigmaz), \quad t\ge 0,
		\end{equation}
		where $x\in\R^n$ is the state, $u_i,\;i = 1, \ldots, m$, are the (scalar) control inputs, $f_i$ and $g_{i, j}$ are twice continuously differentiable vector fields on $\R^n$, with $f_i(0) = 0, g_{i, j}(0) = 0$, for each $i\in\PSet, j\in\{1, \ldots, m\}$. Let $\CSet$ be the set where the control $u := [u_1, \ldots, u_m]\transp$ takes its values. With a feedback control function $\ol k_\sigma(x) = \left[k_{\sigma, 1}(x), \ldots, k_{\sigma, m}(x)\right]\transp$, the closed-loop system stands as:
		\begin{equation}
			\label{e:ssysdefconcl}
			\dot x = f_\sigma(x) + \sum_{j = 1}^m g_{\sigma, j}(x) k_{\sigma, j}(x),\qquad (x(0), \sigma(0)) = (\xz,\sigmaz), \quad t\ge 0.
		\end{equation}
		We now describe the controller design methodology. A universal formula for stabilization of control-affine nonlinear systems was first constructed in~\cite{ref:sontagunivformula}, for the control taking values in $\CSet = \R^m$. The articles~\cite{ref:linunivformulabddcon},~\cite{ref:linunivformularestrcon}, and~\cite{ref:malisoffunivformulaMinkball} provide universal formulae for bounded controls, positive controls, and controls restricted to Minkowski balls, respectively. In view of the analysis results of \secref{s:mainres} and the universal formulae provided in the aforementioned articles, it is possible to synthesize controllers $\ol k_\sigma$ for~\eqref{e:ssysdefcon} such that the closed-loop system~\eqref{e:ssysdefconcl} is \gasas{} In general, we obtain one synthesis scheme for each type of $\CSet$. The following theorem provides a typical illustration of such a result for the case $\CSet = \R^m$; a complete recipe to obtain such results in other cases is provided in Remark~\ref{r:recipe}.

		\begin{theorem}
			Consider the system~\eqref{e:ssysdefcon}, with $\CSet = \R^m$. Suppose that $\sigma$ is of class EH, and there exists a family $\{V_i\}_{i\in\PSet}$ of twice continuously differentiable real-valued functions on $\R^n$, such that
			\begin{enumerate}[label={\rm (C\arabic*)}, align=left, leftmargin=*]
				\item (V1) of Assumption~\ref{a:V} holds;\label{cond:u:1}
				\item (V3) of Assumption~\ref{a:V} holds;\label{cond:u:4}
				\item $\therex\{\lambda_i\}_{i\in\PSet}\subseteq\R$ such that $\fa x\in\R^n\setmin\{0\},\;\fa i\in\PSet$, \label{cond:u:2}
					\[
						\inf_{u\in\CSet}\!\left\{\!\LieD{f_i}{V_i(x)}\! + \!\lambda_i V_i(x) \!+ \!\sum_{j = 1}^m\!\!u_j\LieD{g_{i, j}}{V_i(x)}\right\} < 0;
					\]
				\item $\fa \eps > 0 \;\therex \delta > 0$ such that if $x(\neq 0)$ satisfies $\norm{x} < \delta$, then $\therex u\in\R^m,\; \norm u < \eps$, such that $\fa i\in\PSet$,\footnote{This is known as the small-control property~\cite{ref:sontagunivformula}.}\label{cond:u:3}
					\[
						\LieD{f_i}{V_i} + \sum_{j = 1}^m u_j\cdot\LieD{g_{i, j}}{V_i} \le -\lambda_i V_i;
					\]
				\item (E3)-(E4) of Theorem~\ref{t:gasaspoisson} hold.
			\end{enumerate}
			Then the feedback control function
			\[
				\ol k_\sigma(x) = [k_{\sigma, 1}(x), \ldots, k_{\sigma, m}(x)]\transp,
			\]
			where
			\begin{subequations}
			\begin{align}
				k_{i, j}(x) & := \displaystyle{-\LieD{g_{i, j}}{V_i}(x)\;\varphi\!\left(\ol W_i(x), \wt W_i(x)\right)}\\
				\ol W_i(x) & := \LieD{f_i}{V_i}(x) + \lambda_i V_i(x),\\
				\wt W_i(x) & := \sum_{j = 1}^m\left(\LieD{g_{i, j}}{V_i}(x)\right)^2,\\
			\intertext{and}
				\varphi(a, b) & := 
				\begin{cases}
					\displaystyle{\frac{a + \sqrt{a^2 + b^2}}{b}}\quad & \text{if }b\neq 0,\\
					0 & \text{otherwise,}
				\end{cases}
				\label{e:phidef}
			\end{align}
			\label{e:udefs}
			\end{subequations}
			renders~\eqref{e:ssysdefconcl} \gasas{}
			\label{t:ugasas}
		\end{theorem}
		\begin{proof}
			The proof relies heavily on the construction of the universal formula in~\cite{ref:sontagunivformula}. Fix $t\in\posR$. If $x\neq 0$, applying the definition of $\varphi$, we get
			\begin{align*}
				& \LieD{f_{\sigma(t)}}{V_{\sigma(t)}}(x) + \sum_{i = 1}^m k_{\sigma(t), i}(x)\LieD{g_{\sigma(t), i}}{V_{\sigma(t)}}(x)\\
				& = \LieD{f_{\sigma(t)}}{V_{\sigma(t)}}(x) - \wt W_{\sigma(t)}(x)\!\cdot\!\varphi\!\left(\!\ol W_{\sigma(t)}(x),\!\left(\wt W_{\sigma(t)}(x)\right)^{\!2}\right)\\
				& = -\lambda_{\sigma(t)} V_{\sigma(t)}(x) -\sqrt{\left(\LieD{f_{\sigma(t)}}{V_{\sigma(t)}}(x) + \lambda_{\sigma(t)}V_{\sigma(t)}(x)\right)^2 + \left(\wt W_{\sigma(t)}(x)\right)^2}\\
				& < -\lambda_{\sigma(t)} V_{\sigma(t)}(x).
			\end{align*}
			Since $t$ is arbitrary, we conclude that the above inequality holds for all $t\in\posR$. Note that by (C3), if $x\in\bigcap_{j = 1}^m \ker\left(\LieD{g_{i, j}}{V_i}\right)$ for any $i\in\PSet$, we automatically have $\LieD{f_{\sigma(t)}}{V_{\sigma(t)}}(x) + \lambda_{\sigma(t)} V_{\sigma(t)}(x) < 0$.
			The above arguments, in conjunction with (C1) and (C2) enable us to conclude that the family $\{V_i\}_{i\in\PSet}$ satisfies Assumption~\ref{a:V} for the closed-loop system~\eqref{e:ssysdefconcl}. (E1) and (E2) hold by hypotheses. The assertion now follows from Theorem~\ref{t:gasaspoisson}.
		\end{proof}

		\begin{remark}
			{\rm Theorem~\ref{t:ugasas} can be modified to suit a different control set $\CSet$ and a different type of $\sigma$ using the following simple recipe. First, recall from the discussion preceding Theorem~\ref{t:ugasas} that $\CSet$ may be any one among $\R^m$, the nonnegative orthant of $\R^m$, the unit ball (with respect to the Euclidean norm) of $\R^m$, and a Minkowski ball in $\R^m$. Now suppose that a $\CSet$ is given to us, and let $\sigma$ belong to class UH. Then:
			\begin{enumerate}[label={\rm (R\arabic*)}, align=left, leftmargin=*]
				\item (C1) and (C2) remain unchanged;
				\item the given $\CSet$ replaces the $\CSet = \R^m$ in Theorem~\ref{t:ugasas};
				\item (U3) replaces (E3)-(E4) in (C5);
				\item the universal formula corresponding to the given $\CSet$ replaces the one given in~\eqref{e:udefs}.\RemarkEnd
			\end{enumerate}
			\label{r:recipe}
			}
		\end{remark}

		\subsection{Mode-independent Controllers}
		\label{s:synmi}
		Consider the affine-in-control switched system~\eqref{e:ssysdefcon}. Let $\ol k(x) = \left[k_{1}(x), \ldots, k_{m}(x)\right]\transp$ be a feedback control function, with which the closed-loop system stands as:
		\begin{equation}
			\dot x = f_\sigma(x) + \sum_{j = 1}^m g_{\sigma, j}(x) \ol k_j(x), \qquad (x(0), \sigma(0)) = (\xz, \sigma_0), \quad t\ge 0.
			\label{e:ssysdefconclunob}
		\end{equation}

		\begin{theorem}
			Consider the system~\eqref{e:ssysdefcon} with $\CSet = \R^m$. Suppose that $\sigma$ belongs to class EH, and there exists a family $\{V_i\}_{i\in\PSet}$ of twice continuously differentiable real-valued functions on $\R^n$ such that
			\begin{enumerate}[label={\rm (\roman*)}, align=right, leftmargin=*, widest=iii]
				\item (V1) and (V3) of Assumption~\ref{a:V} holds;
				\item there exists a control function $\ol k:\R^n\lra\CSet$, such that $\LieD{f_i + g_i\ol k}{V_i}(x) \le -\lambda_i V_i(x)$ for every $i\in\PSet$, $x\in\R^n$, for some $\{\lambda_i\}_{i\in\PSet}\subset\R$;
				\item (E3)-(E4) of Theorem~\ref{t:gasaspoisson} holds.
			\end{enumerate}
			Then $\ol k$ renders~\eqref{e:ssysdefcon} \gasas{} in closed-loop.
			\label{t:gasasunobcon}
		\end{theorem}
		
		Note that this result does not need a feedback controller $\ol k$ that simultaneously stabilizes the family~\eqref{e:ssysfam}, which in general is difficult to get; it proposes controllers which may leave some subsystems unstable, but nonetheless achieve \gasas{} of the closed-loop switched system.

		\emph{Proof of Theorem~\ref{t:gasasunobcon}.}
			The assertion follows immediately by first observing that the closed-loop system is~\eqref{e:ssysdefconclunob}, and then applying Theorem~\ref{t:gasaspoisson} to~\eqref{e:ssysdefconclunob}. Indeed, note that hypothesis (ii) holds for~\eqref{e:ssysdefconclunob} by our assumption on $\sigma$, (iii) implies (EH3)-(EH4) hold, and (i)-(ii) ensure that (E1) holds.$\;\;$\endproof

	\section{Conclusion and further work}
	\label{s:concl}
		We have established sufficient conditions for global asymptotic stability almost surely, in the mean, and in probability, of randomly switched systems and a methodology for almost sure global asymptotic stabilization and global asymptotic stabilization in the mean of randomly switched systems with control inputs. The switching signals were assumed to be semi-Markovian.

		An interesting research direction is to extend the above results to systems with disturbance inputs. The analysis becomes more involved, and for synthesis tools universal formulae for \iss{} disturbance attenuation in nonlinear control literature are needed. Some preliminary results have been reported in~\cite{ref:extstab} and~\cite{ref:myphdthesis}. In the particular case of Markovian switching signals, one can prove stochastic analogs of \emph{input to state stability} (\iss{})~\cite[Chapter~3]{ref:myphdthesis}.

	\section*{Acknowledgments}
		The authors are grateful to Sean~P.~Meyn and P.~R.~Kumar for motivating the results of~\secref{s:general}, and V.~S.~Borkar for helpful discussions.


\end{document}